\UseRawInputEncoding

\documentclass[10pt]{amsart}

\newtheorem{theorem}{Theorem}
\newtheorem{lemma}{Lemma}

\newtheorem{cor}{Corollary}
\newtheorem{prop}{Proposition}

\usepackage{amsmath}
\usepackage[psamsfonts]{amssymb}
\usepackage[mathscr]{euscript}

\newcommand{\cal}{\mathcal}

\title[Generalized paracomplex structures
on generalized reflector spaces] {Generalized paracomplex structures
on generalized reflector spaces}
\author{Johann Davidov}

\thanks{The author is partially supported by the National Science Fund,
Ministry of Education and Science of Bulgaria under contract
KP-06-N52/3.}

\address{Institute of Mathematics and Informatics \\
Bulgarian Academy of Sciences\\ Acad. G.Bonchev Str. Bl.8\\
1113 Sofia\\ Bulgaria} \email{jtd@math.bas.bg}

\begin{document}

\begin{abstract}

Non-trivial examples of generalized paracomplex structures (in the
sense of the generalized geometry \`a la Hitchin) are constructed
applying the twistor space construction scheme.

 \vspace{0,1cm} \noindent 2010 {\it Mathematics Subject Classification} 53C28; 53C15, 53D18.

\vspace{0,1cm} \noindent {\it Key words: paracomplex structures,
generalized paracomplex structures, generalized reflector spaces}.

\end{abstract}

\maketitle \vspace{0.5cm}

\section{Introduction}

The seminal work \cite{Hit04} of Nigel Hitchin opened up a new
research area, the generalized complex geometry which has further
been developed by his students M. Gualtieri \cite{Gu}, G. Cavalcanti
\cite{Ca}, F. Witt \cite{Witt}. The concept of a generalized complex
structure encompasses complex, symplectic and holomorphic Poisson
structures. This and the fact that generalized complex structures
are involved in supersymmetric $\sigma$-models motivates the great
interest in generalized complex geometry.

A  generalized almost complex structure on a smooth manifold $M$ is
an endomorphism ${\mathcal J}$ of the bundle $TM\oplus T^{\ast}M$
satisfying ${\mathcal J}^2=-Id$ and compatible with the metric
$<X+\alpha,Y+\beta>=\alpha(Y)+\beta(X)$, the latter metric being of
neutral signature $(n,n)$, $n=dim\,M$.

In line with this definition, A. Wade \cite{Wa} defined the notion
of a generalized paracomplex structure (see also I. Vaisman
\cite{Va07}). Recall first that an almost paracomplex structure on a
smooth manifold $M$ is an endomorphism $K$ of the tangent bundle
$TM$ such that $K^2=id$ and the $\pm 1$-eigen subbundles of $TM$
have the same rank. Such a structure is called integrable if its
Nijenhuis tensor
$$
N(X,Y)=[X,Y]+[KX,KY]-K[KX,Y]-K[X,KY]
$$
vanishes (as is standard,  an integrable almost paracomplex
structure is called paracomplex). An almost paracomplex structure
$K$ is called compatible with a metric $g$ or almost para-Hermitian
if $g(KX,Y)+g(X,KY)=0$. A metric with this property is necessarily
of neutral signature. We refer to \cite{CFG} and \cite{AMT} for a
survey on paracomplex geometry. Now, following \cite{Wa}, define a
generalized almost paracomplex structure as an endomorphism
${\mathcal K}$ of the bundle $TM\oplus T^{\ast}M$ satisfying
${\mathcal K}^2=Id$ and compatible with the metric $<.\,,.>$. The
$\pm 1$-eigen subbundles for such a structure ${\mathcal K}$ are
both of rank $n$. Similar to the case of a usual almost paracomplex
structure, the integrability condition for  a generalized almost
paracomplex structure  is defined as the vanishing of its Nijenhuis
tensor, however this tensor is defined by means of the bracket
introduced by T. Courant \cite{Cou} instead of the Lie bracket. The
generalized paracomplex structures include the paracomplex,
symplectic and Poisson structures as special cases \cite{Wa}.

A specific feature of the generalized geometry is the existence of
the so-called $B$ transform of $TM\oplus T^{\ast}M$ determined by a
skew-symmetric $2$-form on $M$. Such a form $B$ acts on  $TM\oplus
T^{\ast}M$ via the inclusion  $\Lambda^{2}T^{\ast}M\subset
\Lambda^{2}(TM\oplus T^{\ast}M)\cong \frak{so}(<.\,,.>)$ (the space
of endomorphisms skew-symmetric with respect to  the metric
$<.\,,.>$). The resulting map, also denoted by $B$, is $X+\alpha\to
\imath_{X}B$. The exponential map $e^B$ is given by $X+\alpha\to
X+\alpha+\imath_{X}B$ and is an orthogonal transformation of
$TM\oplus T^{\ast}M$ called the  $B$-transform. Thus, if ${\mathcal
K}$ is a generalized complex or paracomplex structure, $e^B{\mathcal
K}e^{-B}$ is also such a structure called the $B$-transform of
${\mathcal K}$. If the form $B$ is closed, ${\mathcal K}$ is
integrable if and only if its $B$-transform is integrable as follows
from  \cite[Proposition 3.23]{Gu}.

The notion of a pseudo-Riemannian metric also has a generalized
version \cite{Gu,Witt, Fer}.  A generalized metric of signature
$(r,s)$ on a  manifold $M$ is a subbundle $E$ of $TM\oplus
T^{\ast}M$ such that $rank\, E=dim\, M$, the restriction of the
metric $<.\,,.>$ to $E$ is of signature $(r,s)$, and $E_x\cap
T_x^{\ast}M=\{0\}$ for every $x\in M$. Every such a bundle is
uniquely determined by a metric $g$ of signature $(r,s)$ and a
skew-symmetric $2$-form $\Theta$ on $M$ such that
$E=\{X+g(X)+\Theta(X):~X\in TM\}$. In the case of signature $(n,0)$,
the positive definite one, the condition $E_x\cap T_x^{\ast}M=\{0\}$
is automatically satisfied.

Given a generalized metric $E=\{X+g(X)+\Theta(X):~X\in TM\}$ of
neutral signature $(n,n)$ on a manifold $M$, denote by ${\mathcal
G}(E)$ the bundle over $M$ whose fibre at a point $p\in M$ consists
of generalized paracomplex structures ${\mathcal K}$ on the tangent
space $T_pM$ compatible with the generalized metric $E_p$, the fibre
of $E$ at $p$, in the sense that ${\mathcal K}E_p\subset E_p$. We
call this bundle generalized reflector space of $(M,E)$. The reason
for this name is that if $(M,g)$ is a pseudo-Riemannian manifold of
neural signature, the bundle over $M$ whose fibres are compatible
paracomplex structures on the tangent spaces of $M$ is called the
reflector space of $(M,g)$. Using the Levi-Civita  connection of the
metric $g$, one can define on the reflector space two almost
paracomplex structures which are analogs of the
Atiyah-Hitchin-Singer \cite{AHS} and Eells-Salamon \cite{ES} almost
complex structures on the twistor space of a Riemannian manifold. In
the case of the generalized reflecter space ${\mathcal G}(E)$, we
use the connection $D'$ on the bundle $E'=E$ defined in Hitchin's
paper \cite{Hit10}. Transferring this connection by means of the
isomorphism $pr_{TM}|E:E\to TM$ gives a connection $\nabla$ on $TM$
compatible with the metric $g$ whose torsion $3$-form is $d\Theta$
[ibid.]. If $E''=E^{\perp}$ is the orthogonal complement of $E$ with
respect to $<.\,,.>$, the projection $pr_{TM}|E'':E''\to TM$ is also
an isomorphism and we can transfer the connection $\nabla$ to a
connection $D''$ on $E''$ by means of this isomorphism. Thus, we
obtain a connection $D$ on the bundle $TM\oplus T^{\ast}M=E'\oplus
E''$ setting $D=D'\oplus D''$. This connection yields  a splitting
${\mathcal V}\oplus {\mathcal H}$ of the tangent bundle of the
manifold ${\mathcal G}(E)$ into vertical and horizontal parts. Then,
following the standard twistor scheme, we can define four
generalized almost paracomplex structures ${\mathscr
K}_{\varepsilon}$ on  ${\mathcal G}(E)$, $\varepsilon=1,\dots,4$.
The structure ${\mathcal K}_1$ is a generalized paracomplex analog
of the Atiyah-Hitchin-Singer almost complex structure, while
${\mathcal K}_2$, ${\mathcal K}_3$, ${\mathcal K}_4$ are analogs of
the Eells-Salamon structure. The construction is similar to that in
\cite{Dav} where the case of a positive definite generalized metric
and the bundle of compatible generalized complex structures is
studied. The differences come mainly from the different signatures
of the generalized metrics under consideration.

In the present paper, it is shown that the structures ${\mathcal
K}_2$, ${\mathcal K}_3$, ${\mathcal K}_4$ are never integrable in
accordance with classical twistor theory. The main result of the
paper gives necessary and sufficient conditions for integrability of
the generalized almost paracomplex structure ${\mathcal K}_1$ in
terms of the curvature of the neutral-signature pseudo-Riemannan
manifold $(M,g)$.

It  seems there are not many examples of generalized (almost)
paracomplex manifolds. The present paper provides non-trivial
examples of  such manifolds, i.e.,  endowed with structures that are
not yielded by a paracomplex, symplectic, or Poisson structure, or a
$B$-transform of a structure of this type.

\medskip

\section{Paracomplex and generalized paracomplex structures}

\subsection{Paracomplex structures on vector spaces}
In this paragraph, we recall basic facts about paracomplex
structures.

Let $T$ be a finite dimensional real vector space. A product
structure on $T$ is an endomorphism $P$ of $T$ with $P^2=Id$,
$P\neq\pm Id$. Denote by $T^{\pm}$ the eigenspaces of $T$ associated
to the eigenvalues $\pm 1$ of $T$. If $dim\,T^{+}=dim\,T^{-}$, then
$P$ is called a paracomplex structure.

A product structure $K$ on $T$ is called compatible with a metric
$g$ on $T$ if it is skew-symmetric with respect to $g$:
$$
g(KX,Y)+g(X,KY)=0,\quad X,Y\in T.
$$
A metric $g$ satisfying this identity is necessarily of neutral
signature $(n,n)$. The pair $(g,K)$ is said to be a para-Hermitian
structure. Given such a structure, if $e_1$ is a vector with
$g(e_1,e_1)=1$, we have $g(Ke_1,e_1)=0$ and $g(Ke_1,Ke_1)=-1$. Then
$W=span\{e_1,Ke_1\}$ is a $2$-dimensional subspace of $T$ invariant
under $K$, and $T=W\oplus W^{\perp}$. Thus, induction  on $dim\,T$
shows 
that $T$ admits a $g$-orthonormal basis of the form
$e_1,...,e_n,e_{n+1}=Ke_1,...,e_{2n}=Ke_n$ with
$g(e_i,e_i)=-g(Ke_i,Ke_i)=1$,  $i=1,...,n$. Set
\begin{equation}\label{a-basis}
a_i=\frac{1}{\sqrt 2}(e_i+Ke_i),\quad a_{n+i}=\frac{1}{\sqrt 2}(e_i-Ke_i),\quad i=1,...,n.
\end{equation}
Then $\{a_i,a_{n+i}\}$ is basis of $T$ such that $a_i\in T^{+}$, $a_{n+i}\in T^{-}$, and
\begin{equation}\label{basis A}
g(a_i,a_{n+j})=\delta_{ij},\quad
g(a_i,a_j)=g(a_{n+i},a_{n+j})=0,\quad i,j=1,...,n.
\end{equation}
In particular, $dim\,T^{+}=dim\,T^{-}$, so $K$ is a paracomplex
structure. Moreover $T^{\pm}$ are totally isotropic subspaces of
dimension $n$ (the maximal one). In other words, $T^{\pm}$ are
Lagrangian subspaces with respect to the symplectic form
$\omega(X,Y)=g(X,KY)$ on $T$. Thus, $(T,\omega, T^{+},T^{-})$ is a
linear bi-Lagrangian structure in the sense of \cite{EST}.

Let $P$ be a paracomplex structure. Take a basis $\{a_i,a_{n+i}\}$
of $T$ such that $a_i\in T^{+}$, $a_{n+i}\in T^{-}$. Then the
paracomplex structure $P$ is compatible with the neutral metric $g$
defined by identities (\ref{basis A}).

\subsection{The space of compatible paracomplex structures on a
vector space}\label{space}

Let $T$ be a real vector space of dimension $2n$ endowed with  a
metric $g$  of signature $(n,n)$. Let $e_1,...,e_{2n}$ be an
orthonormal basis of $T$ such that $g(e_i,e_{n+j})=0$,
$g(e_i,e_j)=\delta_{ij}$, $g(e_{n+i},e_{n+j})=-\delta_{ij}$,
$i,j=1,...,n$. Setting $Ke_i=e_{n+i}$, $Ke_{n+i}=e_i$ defines a
paracomplex structure $K$ on $T$ compatible with metric $g$. The set
$Z(T)$ of all such structures has a natural structure of an embedded
submanifold of the vector space $\frak{so}(g)$ of $g$-skew-symmetric
endomorphisms of $T$. The tangent space of $Z(T)$ at a point $K$
consists of endomorphisms in $\frak{so}(g)$ anti-commuting with $K$.
Then we can define a product structure ${\mathcal K}$ on $Z(T)$
setting
$$
{\mathcal K}V=KV ~ {\rm{for}} ~ V\in T_{K}Z(T).
$$
where, as usual, the juxtaposition of two endomorphisms means
"composition". This product structure is compatible with the
restriction to  $Z(T)$ of the standard metric of $\frak{so}(g)$
$$
G(L_1,L_2)=-\frac{1}{2}Trace\,(L_1L_2).
$$
In particular, ${\mathcal K}$ is a paracomplex structure.

\smallskip

\noindent {\it Notation}. The value of a section $s$ of a vector
bundle at a point $p$ will be denote by $s_p$ or $s(p)$.

\smallskip

Denote by $D$  the Levi-Civita connection of the metric $G$ on
$Z(T)$. Let $X,Y$ be vector fields on $Z(T)$ considered as
$\mathfrak{so}(g)$-valued functions on $\mathfrak{so}(g)$. By the
Koszul formula, for every $K\in Z(T)$,
\begin{equation}\label{coder}
(D_{X}Y)_K=\frac{1}{2}(Y'(K)(X_K)-K\circ Y'(K)(X_K)\circ K)
\end{equation}
where $Y'(K)\in Hom(\mathfrak{so}(g),\mathfrak{so}(g))$ is the
derivative of the function $Y:\mathfrak{so}(g)\to \mathfrak{so}(g)$
at the point $K$. The latter formula implies $(D_X{\mathcal
K})(Y)=0$ since $({\mathcal K}Y)'(K)(X)=X_K\circ Y_K+K\circ
Y'(K)(X_K)$. Thus $(G,{\mathcal K})$ is a para-K\"ahler structure on
$Z(T)$.

\smallskip

The group $O(g)$ of $g$-orthogonal transformations of $T$ acts
transitively on $Z(T)$. The isotropy subgroup at a point $K\in Z(T)$
is isomorphic to the paraunitary group $U(n,{\mathbb A})$, ${\mathbb
A}=\{x+\varepsilon y:~x,y\in{\mathbb R}, \varepsilon^2=1\}$ being
the algebra of paracomplex numbers.

 If $\{a_i,a_{n+i}\}$ is a basis
of $T$ such that
\begin{equation}\label{or}
Ka_i=a_i,\quad Ka_{n+i}=-a_{n+i},\quad g(a_i,a_{n+j})=\delta_{ij},\quad i,j=1,...,n,
\end{equation}
then $g(a_i,a_j)=g(a_{n+i},a_{n+j})=0$ and  every transformation in $O(g)$ commuting with $K$ has the block-matrix
representation
\begin{equation}\label{iso}
\left(\begin{array}{cccc}
    A&0\\
    0&(A^{-1})^t\\
\end{array}\right),\quad A\in GL(n,{\mathbb R})
\end{equation}
with respect to the basis $\{a_i,a_{n+i}\}$. Therefore $Z(T)$ has
the homogeneous representation $Z(T)\cong O(n,n)/U(n,{\mathbb A})$.
In particular, $dim\,Z(T)=n^2-n$.

Given $K\in Z(T)$, the transition matrix of every two bases
satisfying identities (\ref{or}) has the form (\ref{iso}), so they
yield the same orientation on $T$, called the orientation induced by
$K$. Let $e_1,,...,e_{n}, e_{n+1}=Ke_{1}, ...,e_{2n}=Ke_{n}$ be an
orthonormal basis with $g(e_{i},e_{j})=\delta_{ij}$, $i,j=1,...,n$.
Let $\{a_i,a_{n+i}\}$ be the basis of $T$ defined by means of
$\{e_i,e_{n+i}\}$ via (\ref{a-basis}). Then
$$
\begin{array}{c}
a_1\wedge a_2\wedge...\wedge a_{2n-1}\wedge
a_{2n}=(-1)^{\frac{n(n-1)}{2}}a_1\wedge a_{n+1}\wedge\cdots\wedge
a_n\wedge
a_{2n}\\[6pt]
=(-1)^{n} e_1\wedge e_2\wedge\cdots\wedge e_{2n-1}\wedge e_{2n}.
\end{array}
$$
Therefore all bases $\{e_{i},Ke_{i}\}$ with $g(e_i,e_j)=\delta_{ij}$
determine the same orientation which coincides with the orientation
induced by $K$ if $n$ is even.

If we fix an orientation on $T$, the space $Z(T)$ splits into the
spaces
$$
Z_{\pm}(T)=\{K\in Z(T):~ K~\rm{induces}~\pm ~\rm{the ~orientation~
of}~ T\}.
$$
These are the connected components of $Z(T)$ with homogeneous
representation $SO(n,n)/U(n,{\mathbb A})$.

The metric $g$ on $T$ induces a metric on $\Lambda^2 T$ defined by
$$
g(v_1\wedge v_2,v_3\wedge v_4)=g(v_1,v_3)g(v_2,v_4)-g(v_1,v_4)g(v_2,v_3),\quad v_1,...,v_4\in T.
$$
As is standard, using the metric $g$, we identify $\Lambda^2 T$ with the space $\frak{so}(g)$ of $g$-skew-symmetric
endomorphisms of $T$. The endomorphism $S_a$ corresponding to a $2$-vector $a\in\Lambda^2 T$ is defined by
\begin{equation}\label{Sa}
g(S_a u,v)=g(a,u\wedge v),\quad u,v\in T.
\end{equation}

Suppose  $T$ is oriented and of dimension four. Then the Hodge star
operator $\ast$ acts as an involution on $\Lambda^2 T$ and we have
the orthogonal decomposition $\Lambda^2T=\Lambda^2_{+}T\oplus
\Lambda^2_{-}T$ into the subspaces corresponding to the eigenvalues
$\pm 1$ of the Hodge operator. Let $(e_1,...,e_4)$ be an oriented
orthonormal basis of $T$ with $g(e_1,e_1)=g(e_2,e_2)=1$,
$g(e_3,e_3)=g(e_4,e_4)=-1$. Then
$$
\ast({\bf e}_1\wedge{\bf e}_2)={\bf e}_3\wedge{\bf e}_4,\quad
\ast({\bf e}_1\wedge{\bf e}_3)={\bf e}_2\wedge{\bf e}_4,\quad
\ast({\bf e}_1\wedge{\bf e}_4)=-{\bf e}_2\wedge{\bf e}_3.
$$
Set
\begin{equation}\label{s-basis}
s_1^{\pm}=\frac{1}{\sqrt 2}({\bf e}_1\wedge{\bf e}_2\pm{\bf e}_3\wedge{\bf e}_4),\quad s_2^{\pm}=\frac{1}{\sqrt 2}({\bf
e}_1\wedge{\bf e}_3\mp {\bf e}_4\wedge{\bf e}_2),\quad s_3^{\pm}=\frac{1}{\sqrt 2}({\bf e}_1\wedge{\bf e}_4\mp{\bf
e}_2\wedge{\bf e}_3).
\end{equation}
The bi-vectors $\{s_1^\pm,s_2^\pm,s_3^\pm\}$ constitute an
orthogonal basis of $\Lambda^2_{\pm}T$ and
$||s_1^{\pm}||^2=1$,$||s_2^{\pm}||^2=||s_3^{\pm}||^2=-1$. The
endomorphisms $J_i=S_{s_i^{+}}$ corresponding to ${\sqrt 2}s_i^{+}$
via (\ref{Sa}), $i=1,2,3$, are skew-symmetric, anti-commute with
each other and satisfy the identities $J_1^2=-Id$, $J_2^2=J^2_3=Id$,
$J_3=J_2J_1$. Hence, for $a\in\Lambda^2_{+}T$,
$K_a=y_1J_1+y_2J_2+y_3J_3$ is a skew-symmetric endomorphism of $T$
which is a product structure if and only if $-y_1^2+y_2^2+y_3^2=1$.
If the latter identity is satisfied, $K_a$ is an ant-symmetric
product structure, i.e., a paracomplex structure.  Suppose
$-y_1^2+y_2^2+y_3^2=1$, and set
$$
u_1=e_1,\quad
u_2=\frac{1}{\sqrt{y_2^2+y_3^2}}[(y_2^2+y_3^2)e_2+y_1y_2e_3+y_1y_3e_4],\quad
u_3=K_au_1,\quad u_4=K_au_2.
$$
Then  $u_1,u_2,u_3=K_au_1,u_4=K_au_2$ is an orthonormal basis of $T$
with $g(u_1,u_1)=g(u_2,u_2)=1$,  and $u_3=y_1e_2+y_2e_3+y_3e_4$,
$u_4=\frac{1}{\sqrt{y_2^2+y_3^2}}(-y_3e_3+y_2e_4)$ . This basis
yields the  orientation of $T$,  $u_1\wedge...\wedge
u_4=e_1\wedge...\wedge e_4$, hence $K_a$ induces the orientation of
$T$. Thus, $Z_{+}(T)$ is the one-sheeted hyperboloid
$-y_1^2+y_2^2+y_3^2=1$ in the $3$-dimensional vector space
$\Lambda^2_{+}T$.

\medskip

\noindent {\it Remark}. Similarly, one can see that the space of
complex structures on $T$ compatible with the matric $g$ and
inducing the orientation of $T$ is the two-sheeted hyperboloid
$-y_1^2+y_2^2+y_3^2=-1$ in $\Lambda^2_{+}T$ (\cite{BDM-05}).

\subsection{The reflector space of a manifold with a neutral
metric}\label{reflector}

Recall first that an almost product structure $P$ on a manifold $N$
is an endomorphism $P$ of the tangent bundle $TN$ such that
$P^2=Id$, $P\neq\pm Id$ at every point. An almost product structure
is called integrable  if its Nijenhuis tensor
$$
N(X,Y)=[X,Y]+[PX,PY]-P[PX,Y]-P[X,PY]
$$
vanishes. An almost product structure $P$ is called almost
paracomplex if  the subbundles of $TN$ associated to the eigenvalues
$\pm 1$ of $P$ have the same rank.  We refer to \cite{AMT, CFG} for
a survey on paracomplex geometry.

Let $(M,g)$ be a pseudo-Riemannian manifold with a metric $g$ of
neutral signature $(n,n)$. Denote by ${\mathcal Z}$ the bundle over
$M$ whose fibre at a point $p\in M$ consists of paracomplex
structures on the tangent space $T_pM$ compatible with the metric
$g_p$. This bundle is an analog of the twistor bundle of compatible
complex structures on the tangents space of a (pseudo) Riemannian
manifold. It is called the reflector space of $(M,g)$ in \cite{JR}.
Using the standard scheme of the twistor theory, one can define two
almost paracomplex structures ${\mathcal K}^1$ and ${\mathcal K}^2$
which are analogs of the Atiyah-Hitchin-Singer \cite{AHS} and
Eells-Salamon \cite{ES} almost complex structures on the twistor
space of a Riemannian manifold. To recall the definition of these
almost paracomplex structure, consider ${\mathcal Z}$ as a subbundle
of the vector  bundle $A(TM)$ of $g$-skew-symmetric endomorphisms of
$TM$ considered with the connection induced by the Levi-Civita
connection of $g$. Then, as it is easy to see, the horizontal bundle
of $A(TM)$ is tangent to ${\mathcal Z}$. Thus, we have the
decomposition $T{\mathcal Z}={\mathcal V}\oplus {\mathcal H}$ where
${\mathcal V}$ is the vertical bundle of ${\mathcal Z}$. For every
$K\in{\mathcal Z}$, the vertical space ${\mathcal V}_K$ is tangent
to the fibre of ${\mathcal Z}$ through $K$. The latter has a
paracomplex structure ${\mathcal K}$ described in Sec.~\ref{space}
and we set ${\mathcal K}^iV=(-1)^{i+1}{\mathcal K}V$ for $V\in
{\mathcal V}_K$, so ${\mathcal K}^iV=(-1)^{i+1}KV$, $i=1,2$. If
$X^h_K\in {\mathcal H}_K$ is the horizontal lift of a vector $X\in
TM$, we set ${\mathcal K}^iX^h_K=(KX)^h_K$.

If $M$ is oriented, denote by ${\mathcal Z}_{\pm}\to M$ the bundle
whose fibre at $p\in M$ is the space of paracomplex structures on
$T_pM$ compatible with the metric and $\pm$ the orientation of
$T_pM$. Following the usual terminology of the twistor theory, we
call the bundles ${\mathcal Z}_{+}$ and ${\mathcal Z}_{-}$ positive
and, respectively, negative reflector space. If $M$ is connected,
these are the connected components of ${\mathcal Z}$. It is proved
in \cite{JR} that if $dim\,M=4$, the restriction of the almost
paracomplex structure ${\mathcal K}^1$ to ${\mathcal Z}_{+}$ is
integrable if and only if $(M,g)$ is anti-self-dual, while
${\mathcal K}^2|{\mathcal Z}_{+}$ is never integrable. We refer to
\cite{JR} for more geometric properties of the reflector spaces over
four-manifolds.

\medskip

\noindent {\it Remark}. If $M$ is oriented, one can also consider the bundles over $M$ whose fibre at a point $p\in M$ is
the space of complex structures on $T_pM$ compatible with metric and $\pm$ the orientation of $T_pM$. They admit two
natural almost complex structures. If $dim\,M=4$ the fibre is a two-sheeted hyperboloid, so each of these bundles is
called in \cite{BDM-05} the hyperbolic twistor space of $(M,g)$ (positive and negative, respectively). We refer to
\cite{BDM-04} and \cite{BDM-05} for geometry of hyperbolic twistor spaces.

\subsection{Generalized paracomplex structures on a vector space} Recall that one of
the basic {\bf formal} principle of generalized geometry at the
vector space level is to replace a vector space $T$ with the space
$T\oplus T^{\ast}$ endowed with the neutral metric
$$
<X+\alpha,Y+\beta>=\frac{1}{2}[\alpha(Y)+\beta(X)],\quad X,Y\in T, \alpha,\beta\in T^{\ast}.
$$
By this principle, a generalized paracomplex structure on $T$ is
defined as a product structure $K$  of $T\oplus T^{\ast}$ compatible
with the metric $<.\,,.>$: $<KA,B>+<A,KB>=0$ for $A,B\in T\oplus
T^{\ast}$.

\smallskip

\noindent {\bf Examples} \cite{Wa}.

\smallskip

\noindent {\bf 1}. $K(X+\alpha)=X-\alpha$ is the so-called trivial
generalized paracomplex structure.

\smallskip

\noindent {\bf 2}. If $\omega$ is a symplectic form on $T$, i.e.,  a
non-degenerate skew-symmetric $2$-form, set
$K_{\omega}(X+\alpha)=\omega^{-1}(\alpha)+\omega(X)$, where
$\omega:T\to T^{\ast}$ is the isomorphism determined by the form
$\omega$.

\smallskip

\noindent {\bf 3}. If $\pi$ is a $2$-vector, set
$K_{\pi}(X+\alpha)=(X-\imath_{\alpha}\pi)-\alpha$, where
$\imath_{\alpha}\pi\in T$ is determined by the identity
$\beta(\imath_{\alpha}\pi)=(\alpha\wedge\beta)(\pi)$ for every
$\beta\in T^{\ast}$.
\smallskip

\noindent 4. If $P$ is a product structure on $T$,
$K_{P}(X+\alpha)=PX-P^{\ast}\alpha$ is a generalized paracomplex
structure.

\smallskip

\noindent {\it Remark}. Examples 1 and 3 have no analog in
generalized complex geometry. As for Example 3, recall that one can
define a generalized complex structure by means of a $2$-vector if,
in addition, a complex structure is given. Note also that if a
vector space admits a generalized complex structure it is of even
dimension \cite{Gu}. Examples 1 and 4 show that this is not true in
the case of generalized paracomplex structure.

Recall that any symplectic form $\omega$ determines a $2$-vector
$\pi$ defined by means of the identity
$(\alpha\wedge\beta)(\pi)=\omega(\omega^{-1}(\alpha),\omega^{-1}(\beta))$
for every $\alpha,\beta\in T^{\ast}$. For such a $2$-vector, the
generalized paracomplex structure of Example 3 takes the form
$K_{\pi}(X+\alpha)=X+\omega^{-1}(\alpha)-\alpha$ (to be compared
with Example 2).

\smallskip

Every automorphism $f$ of $T$ gives rise to an automorphism
$F=f\oplus (f^{-1})^{\ast}$ of $T\oplus T^{\ast}$ which is
orthogonal with respect to the metric $<.\,,.>$. Thus, if $K$ is a
generalized paracomplex structure, $FKF^{-1}$ is so. As in the case
of generalized complex structures, the generalized paracomplex
structures admit another type of symmetries, the $B$-transforms.
Recall that  any $2$-form $B\in \Lambda^{2}T^{\ast}$ acts on
$T\oplus T^{\ast}$ via the inclusion $\Lambda^{2}T^{\ast}\subset
\Lambda^{2}(T\oplus T^{\ast})\cong \frak{so}(<.\,,.>)$; this is the
map $X+\alpha\to \imath_{X}B$. Denote the latter map again by $B$.
Then the invertible map $e^{B}$ is given by $X+\alpha\to
X+\alpha+\imath_{X}B$ and is an orthogonal transformation of
$T\oplus T^{\ast}$ called the {\it $B$-transform}. Thus, given a
generalized paracomplex structure $K$ on $T$, the map $e^{B}Ke^{-B}$
is also a generalized complex structure on $T$, called {\it the
$B$-transform of $K$}.

\subsection {Generalized metrics on a vector space}

\noindent\\
{\bf Definition}. A generalized metric of signature $(r,s)$ on a vector space $T$ is a subspace $E$ of $T\oplus T^{\ast}$
such that

\noindent $(1)$~$dim\,E=dim\,T$;

\noindent $(2)$~The restriction of the metric
$$
<X+\alpha,Y+\beta>=\frac{1}{2}[\alpha(Y)+\beta(X)]
$$
to $E$ is of signature $(r,s)$;

\noindent $(3)$~ $E\cap T^{\ast}=\{0\}$.

The third condition is automatically satisfied in the case of
signature $(n,0)$ (positive definite metric). In this case, the
notion of a generalized metric has been introduced in
\cite{Gu,Witt}. The definition presented here is a particular case
of the general notion of an admissible generalized metric on a
Courant algebroid given in \cite{Fer}.

Conditions $(1)$ and $(3)$ imply that $T\oplus T^{\ast}=E\oplus T^{\ast}$. It follows that  $E$ is the graph of a map
$T\to T^{\ast}$. Let $g$ and $\Theta$ be the bilinear forms on $M$ determined by the symmetric and skew-symmetric parts
of this map.

\smallskip

\noindent {\it Notation}. The map $T\to T^{\ast}$ determined by a bilinear form $\varphi$ on $T$ will be denoted again by
$\varphi$; thus $\varphi(X)(Y)=\varphi(X,Y)$.

\smallskip

Then
$$
E=\{X+g(X)+\Theta(X):~X\in T\}.
$$
The restriction of the metric $<.\,, .>$ to $E$ is
\begin{equation}\label{pd}
<X+g(X)+\Theta(X),Y+g(Y)+\Theta(Y)>=g(X,Y),\quad X,Y\in T.
\end{equation}
Thus $g$ is metric on $M$ of signature $(r,s)$.

Set $E'=E$ and $E''=E^{\perp}$, the orthogonal complement being taken with respect to the metric $<.\,,.>$. It easy to
check that
$$
E''=\{X-g(X)+\Theta(X):~X\in T\}.
$$
Then we have the orthogonal decomposition $T\oplus T^{\ast}=E'\oplus E''$. Moreover it is obvious that $E''\cap
T^{\ast}=\{0\}$.

Note that the $E'$ and $E''$-components of a vector $X\in T$ are
\begin{equation}\label{prX}
\begin{array}{c}
X_{E'}=\displaystyle{\frac{1}{2}}\{X- (g^{-1}\circ \Theta)(X)+ g(X)-
(\Theta\circ g^{-1}\circ \Theta)(X)\},\\[8pt]
X_{E''}=\displaystyle{\frac{1}{2}}\{X+ (g^{-1}\circ \Theta)(X)- g(X)+ (\Theta\circ g^{-1}\circ \Theta)(X)\}
\end{array}
\end{equation}
and the  components  of a $1$-form $\alpha\in T^{\ast}$ are given by
\begin{equation}\label{pralpha}
\begin{array}{c}
\alpha_{E'}=\displaystyle{\frac{1}{2}}\{g^{-1}(\alpha)+\alpha+(\Theta\circ
g^{-1})(\alpha)\},\\[8pt]
\alpha_{E''}=\displaystyle{\frac{1}{2}}\{-g^{-1}(\alpha)+\alpha-(\Theta\circ g^{-1})(\alpha)\}.
\end{array}
\end{equation}

Let $pr_{T}: T\oplus T^{\ast}\to T$ be the natural projection. According to condition $(3)$, the restriction of this
projection to $E$ is an isomorphism. Identity (\ref{pd}) tells us that this isomorphism is an isometry when $E$ is
equipped with the metric $<.\,,.>|E$ and $T$ with the metric $g$. Similarly, the map $pr_{T}|E''$ is an isometry of the
metrics $<.\,,.>|E''$ and $-g$.

\subsection{Generalized paracomplex structures compatible with a generalized metric on a vector space}

Let $T$ is a real vector space of dimension $m=2n$ and let $E$ be a generalized metric on $T$ of signature $(n,n)$.
Denote by $g$ and $\Theta$ the $(n,n)$-signature metric and the skew-symmetric $2$-form on $T$ determined by $E$.

As in \cite{Dav}, we shall say that a generalized paracomplex
structure $K$ on $T$ is {\it compatible} with the generalized metric
$E$ if $K$ preserves the space $E$.  In fact, this definition is
implicitly given in \cite{Gu} in the case of a generalized complex
structure, see also \cite{Ca-book}.

If $K$ is a such structure, $KE'=E'$ and $KE''=E''$, where, as
above, $E'=E$ and $E''=E^{\perp}$. Setting
\begin{equation}\label{K1,2}
K_1=(pr_{T}|E')\circ K\circ (pr_{T}|E')^{-1},\quad K_2=(pr_{T}|E'')\circ K\circ (pr_{T}|E'')^{-1}.
\end{equation}
we obtain two paracomplex structures on $T$ compatible with neutral metric $g$ determined by $E$. Then
\begin{equation}\label{P12-P}
\begin{array}{c}
K(X+g(X)+\Theta(X))=K_1X+g(K_1X)+\Theta(K_1X),\\[6pt]
K(X-g(X)+\Theta(X))=K_2X-g(K_2X)+\Theta(K_2X).
\end{array}
\end{equation}
Thus, similar to the case of a generalized complex structure (\cite{Gu}), the generalized paracomplex structure $K$ can
be reconstructed from the data $g,\Theta,K_1,K_2$.

\begin{prop}\label{K-K1,2}
Let $g$ be a metric of neutral signature, $\Theta$ - a skew-symmetric $2$-form on $T$, and $K_1$, $K_2$ - two paracomplex
structures compatible with the metric $g$. Let $\omega_{s}(X,Y)=g(X,K_{s}Y)$ be the fundamental $2$-forms of the
para-Hermitian structure $(g,K_s)$, $s=1,2$. Then the block-matrix representation of the compatible generalized
paracomplex structure $K$ determined by  the data $(g,\Theta,K_1,K_2)$ is of the form
$$
K=\frac{1}{2}\left(
  \begin{array}{cc}
    I & 0 \\
    \Theta & I \\
  \end{array}
  \right)
  \left(
  \begin{array}{cc}
    K_1+K_2 & -\omega_{1}^{-1}+\omega_{2}^{-1} \\
    -\omega_{1}+\omega_{2}& -K_1^{\ast}-K_2^{\ast}) \\
  \end{array}
  \right)
    \left(
  \begin{array}{cc}
    I & 0 \\
    -\Theta & I \\
  \end{array}
  \right),
  $$
where $I$ is the identity matrix and $\Theta$, $\omega_1$, $\omega_2$ stand for the maps $T\to T^{\ast}$ determined by
the corresponding $2$-forms.
\end{prop}
This is an analog of \cite[identity (6.3)]{Gu} and can be proved by
means of identities (\ref{prX}) and (\ref{pralpha}), and the
identities   $K_s\circ g^{-1}=-\omega_s^{-1}$, $g\circ
K_s=-\omega_s$, $g\circ K_s\circ g^{-1}=-K^{\ast}_s$, $s=1,2$.

\smallskip

\noindent {\bf Examples}.

\smallskip

\noindent {\bf 3}. The trivial generalized paracomplex structure
$K(X+\alpha)=X-\alpha$ is not compatible with any generalized
metric.

\smallskip

\noindent {\bf 4}. Let $K_{\omega}$ be the generalized paracomplex
structure defined by a non-degenerate skew-symmetric $2$-form
$\omega$: $K_{\omega}(X+\alpha)=\omega^{-1}(\alpha)+\omega(X)$. Take
a Darboux basis $\{a_i,a_{n+i}\}$ for $\omega$:
$\omega(a_i,a_j)=\omega(a_{n+i},a_{n+j})=0$,
$\omega(a_i,a_{n+j})=\delta_{ij}$. Let $g$ be the metric on $T$ for
which $g(a_i,a_j)=g(a_{n+i},a_{n+j})=0$,
$g(a_i,a_{n+j})=g(a_{n+j},a_i)=\delta_{ij}$. This is a neutral
metric, and the generalized paracomplex structure $K_{\omega}$ is
compatible with the generalized metric $E=\{X+g(X):~X\in T\}$.

\smallskip

Suppose  $K_{\omega}$ is compatible with a generalized metric
$E=\{X+g(X)+\Theta(X):~ X\in T\}$. Let $K_1$ be the paracomplex
structures defined by (\ref{K1,2}) with $K=K_{\omega}$; this
structure is compatible with the metric $g$. Then
$K_1X=\omega^{-1}(g(X)+\Theta(X))$. This identity is equivalent to
$g(X,Y)=\frac{1}{2}(\omega(K_1X,Y)-\omega(X,K_1Y))$,
$\Theta(X,Y)=\frac{1}{2}(\omega(K_1X,Y)+\omega(X,K_1Y))$.
Conversely, suppose  there is a product structure $L$ on $T$ such
that the symmetric bilinear form
$g(X,Y)=\frac{1}{2}(\omega(LX,Y)-\omega(X,LY))$ is non-degenerate.
Then $L$ is compatible with the metric $g$, hence $g$ is of neutral
signature, thus $(g,L)$ is a para-Hermitian structure. Setting
$\Theta(X,Y)=\frac{1}{2}(\omega(LX,Y)+\omega(X,LY))$ defines a
skew-symmetric bilinear form such that $K_{\omega}$ is compatible
with generalized metric $E=\{X+g(X)+\Theta(X):~ X\in T\}$.

A generalized paracomplex structure  preserves a generalized metric
$E=\{X+g(X)+\Theta(X):~ X\in T\}$ if and only if it preserves
$E^{\perp}=\{X-g(X)+\Theta(X):~ X\in T\}$. Hence the generalized
complex structure $K_{\omega}$ preserves $E$ if and only if there
exists a product structure $L$ such that
$g(X,Y)=-\frac{1}{2}(\omega(LX,Y)-\omega(X,LY))$,
$\Theta(X,Y)=\frac{1}{2}(\omega(LX,Y)+\omega(X,LY))$.

\smallskip

\noindent {\bf 5}. Consider the generalized paracomplex structure
$K_{\pi}(X+\alpha)=X-\imath_{\alpha}\pi-\alpha$, where  $\pi\in
\Lambda^2T$.

If $\pi=0$,  $K_{\pi}$ is the trivial generalized paracomplex
structure which is not compatible with any generalized metric.

 Take a para-Hermitian structure $(g,L)$. Let $g^{\wedge}$ be
the metric on $\Lambda^2T$ determined by $g$, and define a
$2$-vector $\pi$ via the identity $g^{\wedge}(\pi,X\wedge
Y)=g(X,LY)$. Then $\imath_{\alpha}\pi=-Lg^{-1}(\alpha)$. It is easy
to see that, if $\Theta(X,Y)=g(X,LY)$, the generalized paracomplex
structure $K_{\pi}$ is compatible with the metric
$E=\{X+g(X)+\Theta(X)\}$.

For simplicity, suppose $dim\,T=4$ and $\pi=e_1\wedge e_2\neq 0$ for
some vectors $e_1,e_2\in T$. The generalized paracomplex structure
$K_{\pi}$ is compatible with a generalized metric
$E=\{X+g(X)+\Theta(X)\}$ if and only if for every $X,Y\in T$
$$
\begin{array}{l}
2g(X,Y)+g(X,e_1)\Theta(Y,e_2)+g(Y,e_1)\Theta(X,e_2)
-\Theta(X,e_1)g(Y,e_2)-\Theta(Y,e_1)g(X,e_2)\\[8pt]

+2\Theta(X,Y)-g(X,e_1)g(Y,e_2)+g(Y,e_1)g(X,e_2)
+\Theta(X,e_1)\Theta(Y,e_2)-\Theta(Y,e_1)\Theta(X,e_2)\\[6pt]

=0.
\end{array}
$$
The left-hand side of this identity is the sum of a symmetric
bilinear form and a skew-symmetric one. Hence it is equivalent to
the following identities
\begin{equation}\label{eq-sym}
2g(X,Y)+g(X,e_1)\Theta(Y,e_2)+g(Y,e_1)\Theta(X,e_2)
-\Theta(X,e_1)g(Y,e_2)-\Theta(Y,e_1)g(X,e_2)=0.
\end{equation}
\begin{equation}\label{eq-skew}
2\Theta(X,Y)-g(X,e_1)g(Y,e_2)+g(Y,e_1)g(X,e_2)
+\Theta(X,e_1)\Theta(Y,e_2)-\Theta(Y,e_1)\Theta(X,e_2)=0
\end{equation}
Identity (\ref{eq-sym}) for $(X,Y)=(e_i,e_j)$, $i,j=1,2$,  reads as
$$
g(e_i,e_j)[1+\Theta(e_1,e_2)]=0.
$$
Identity (\ref{eq-skew}) gives
$$
\Theta(e_1,e_2)[2+\Theta(e_1,e_2)]=0.
$$
If $1+\Theta(e_1,e_2)=0$, then $\Theta(e_1,e_2)=0$ by the latter
identity. Hence $1+\Theta(e_1,e_2)\neq 0$, thus $g(e_i,e_j)=0$. Then
if $\Theta(e_1,e_2)=0$, it follows from (\ref{eq-sym}) that
$g(e_1,Y)=0$ for every $Y\in T$, therefore $e_1=0$ and $\pi=0$, a
contradiction. Hence $\Theta(e_1,e_2)\neq 0$, thus
$2+\Theta(e_1,e_2)=0$. Now, the two-dimensional space
$F=span\{e_1,e_2\}$ is maximal isotropic and there exists a maximal
isotropic subspace $F'$ such that $T=F\oplus F'$. The metric $g$
yields an isomorphism $g:F'\cong F^{\ast}$ defined by
$g(Z)(X)=g(Z,X)$. If $e_1^{\ast},e_2^{\ast}$ is the basis of
$F^{\ast}$ dual to $e_1,e_2$, set $f_1=g^{-1}(e_1^{\ast})$,
$f_2=g^{-1}(e_2^{\ast})$. In this way we get a basis
$e_1,e_2,f_1,f_2$ of $T$ such that $g(e_i,e_j)=g(f_i,f_j)=0$,
$g(e_i,f_j)=\delta_{ij}$. A straightforward verification shows that
identities (\ref{eq-sym}) and (\ref{eq-skew}) are satisfied for
$(X,Y)=(e_i,f_j)$, $i,j=1,2$. Furthermore, these identities are
satisfied for $(X,Y)=(f_i,f_j)$ exactly when
$\Theta(e_1,f_2)=\Theta(e_2,f_1)=0$,
$\Theta(e_1,f_1)=\Theta(e_2,f_2)$,
$2\Theta(f_1,f_2)=\Theta(e_1,f_1)\Theta(e_2,f_2)+1$.

Now, given $\pi=e_1\wedge e_2$, complete $e_1,e_2$ to a basis
$e_1,e_2,f_1,f_2$ of $T$. Define a neutral metric $g$ on $T$  by
$g(e_i,e_j)=g(f_i,f_j)=0$, $g(e_i,f_j)=\delta_{ij}$, and take any
skew-symmetric $2$-form $\Theta$ for which $\Theta(e_1,e_2)=-2$,
$\Theta(e_1,f_2)=\Theta(e_2,f_1)=0$,
$\Theta(e_1,f_1)=\Theta(e_2,f_2)$,
$2\Theta(f_1,f_2)=\Theta(e_1,f_1)\Theta(e_2,f_2)+1$. Then, by the
considerations above, the generalized paracomplex structure
$K_{\pi}$ is compatible with the generalized metric
$E=\{X+g(X)+\Theta(X)\}$

\smallskip

\noindent {\bf 6}. Let $P$ be a paracomplex structure compatible
with a metric $g$. Then the generalized paracomplex structure
$K_P(X+\alpha)=PX-P^{\ast}\alpha$ is compatible with a generalized
metric $E=\{X+g(X)+\Theta(X)\}$ if and only if
$\Theta(PX,Y)+\Theta(X,PY)=0$. Clearly, the form
$\Theta(X,Y)=g(X,PY)$ satisfies the latter condition.

\smallskip

If $g$ is a metric on $T$, its graph $E=\{X+g(X):~X\in T\}$ is a
generalized metric of the same signature. As for a generalized
complex structure (\cite{Dav}), we have the following.

\begin{prop}
Let $g$ be  a metric of signature $(n,n)$ on $T$ and $g^{\ast}$ the metric on $T^{\ast}$ determined by $g$. A generalized
paracomplex structure $K$ on $T$ is compatible with the generalized metric $E=\{X+g(X):~ X\in T\}$ if and only it is
compatible with the neutral metric $\widehat{g}=g\oplus g^{\ast}$ on $T\oplus T^{\ast}$.
\end{prop}

\begin{proof}
Let ${\mathscr G}$ be the endomorphism of $T\oplus T^{\ast}$ defined by
$$
<{\mathscr G}(A),B>=\frac{1}{2}\widehat{g}(A,B),\quad A,B\in T\oplus T^{\ast}.
$$
In fact, ${\mathscr G}$ is given by ${\mathscr
G}(X+\alpha)=g^{-1}(\alpha)+X$, $X\in T$, $\alpha\in T^{\ast}$.
Hence a vector $A$ lies in $E$ if and only if ${\mathscr G}(A)=A$;
similarly $A\in E''=E^{\perp}$ if and only if ${\mathscr G}(A)=-A$.
Moreover, $\widehat{g}({\mathscr G}(A),B)=\widehat{g}(A,{\mathscr
G}(B))$.

Suppose  $K$ is compatible with the metric $\widehat{g}$. Then
$$
\begin{array}{c}
\widehat{g}({\mathscr G}(KA),B)=\widehat{g}(KA,{\mathscr G}(B))=-\widehat{g}(A,K{\mathscr G}(B)\\[6pt]
\widehat{g}(K{\mathscr G}(A),B)=-\widehat{g}({\mathscr G}(A),KB)=-\widehat{g}(A,{\mathscr G}(KB)).
\end{array}
$$
Therefor $K{\mathscr G}(B)={\mathscr G}(KB)$ for every $B\in T\oplus T^{\ast}$. Thus if $B\in E$, then $KB$ also lies in
$E$.

Conversely, suppose that $KE\subset E$. Then $KE^{\perp}\subset
E^{\perp}$, where $E^{\perp}=\{X-g(X):~X\in T\}$ and $T\oplus
T=E\oplus E^{\perp}$. Note that $<.\,,.>|E=2\widehat{g}$,~
$<.\,,.>|E^{\perp}=-2\widehat{g}$ and $\widehat{g}(E,E^{\perp})=0$.

Hence, if $A,B\in E$ or $A,B\in E^{\perp}$,
$$
2[\widehat{g}(KA,B)+\widehat{g}(A,KB)]=\pm[<KA,B>+<A,KB>]=0.
$$
If $A\in E$ and $B\in E^{\perp}$, $\widehat{g}(KA,B)+\widehat{g}(A,KB)]=0-0=0$.

Therefore $K$ is compatible with $\widehat{g}$.
\end{proof}

\subsection{The space of compatible generalized paracomplex structures on a vector
space}\label{CGPCS}

Let $E$ be a generalized metric of signature  $(n,n)$ on a vector
space $T$, $dim\,T=2n$. As above, set $E'=E$ and $E''=E^{\perp}$,
the orthogonal complement being taken with respect to the metric
$<.\,,.>$ on $T\oplus T^{\ast}$.

Denote by ${\mathcal G}(E)$ the set of generalized paracomplex
structures compatible with $E$. This set (non-empty by
Proposition~\ref{K-K1,2}) has the structure of an embedded
submanifold of the vector space $\frak{so}(<.\,,.>)$ of the
endomorphisms of $T\oplus T^{\ast}$, which are skew-symmetric with
respect to the metric $<.\,,.>$. The tangent space of ${\mathcal
G}(E)$ at a point $K$ consists of the endomorphisms $V$ of $T\oplus
T^{\ast}$ anti-commuting with $K$, skew-symmetric with respect to
$<.\,,.>$, and such that $VE\subset E$. Such an endomorphism $V$
sends  $E''$ into itself. Note also that the smooth manifold
${\mathcal G}(E)$ admits a natural paracomplex structure ${\mathscr
P}$  given by $V\to KV$, compatible with the restriction to
${\mathcal G}(E)$ of the standard metric of $\frak{so}(<.\,,.>)$.

For every $K\in {\mathcal G}(E)$, the restrictions $K'=K|E'$ and
$K''=K|E''$ are paracomplex structures on the vector spaces $E'$ and
$E''$ compatible with the neutral metrics $g'= <.\,,.>|E'$ and
$g''=-<.\,,.>|E''$, respectively. Denote by $Z(E')$ and $Z(E'')$ the
sets of paracomplex structures on $E'$ and $E''$ compatible with the
metrics $g'$ and $g''$. Consider these sets with their natural
structures of embedded submanifolds of the  vector spaces
$\frak{so}(g')$ and $\frak{so}(g'')$. The tangent space of $Z(E')$
at $K'$ is $T_{K'}Z(E')=\{V'\in \frak{so}(g'):~V'K'+K'V'=0\}$;
similarly for the tangent space $T_{K''}Z(E'')$. Recall that the
manifold $Z(E')$ admits a paracomplex structure ${\mathscr K}'$
defined by $V'\to K'V'$; similarly $V''\to K''V''$ is a paracomplex
structure ${\mathscr K}''$ on $Z(E'')$. The map $K\to (K',K'')$ is a
diffeomorphism of ${\mathcal G}(E)$  onto $Z(E')\times Z(E'')$
sending a tangent vector $V$ at $K$ to the tangent vector
$(V',V'')$, where $V'=V|E'$ and $V''=V|E''$. Thus ${\mathcal
G}(E)\cong Z(E')\times Z(E'')$ admits four paracomplex structure
${\mathscr P}_{\varepsilon}$ defined by
$$
\begin{array}{c}
{\mathscr P}_1(V',V'')=(K'V',K''V''),\quad {\mathscr P}_2(V',V'')=(K'V',-K''V'')\\[8pt]
{\mathscr P}_3=-{\mathscr P}_2,\quad {\mathscr P}_4=-{\mathscr P}_1.
\end{array}
$$
Clearly, the map $K\to (K',K'')$ is bi-paraholomorphic with respect
to the paracomplex structures ${\mathscr P}$ on ${\mathcal G}(E)$
and ${\mathscr P}_1$ on $Z(E')\times Z(E'')$.

Let $G'(S_1',S_2')=-\frac{1}{2}Trace\,(S_1'S_2')$ be the standard
metric on $\frak{so}(g')$ induced by $g'$; similarly denote by $G''$
the metric on $\frak{so}(g'')$ induced by $g''$. Then, as we have
observed, $(G',{\mathscr K}')$ and $(G'',{\mathscr K}'')$ are
para-K\"ahler structures on $Z(E')$ and $Z(E'')$, so
$(G=G'+G'',{\mathscr P}_{\varepsilon})$, $\varepsilon=1,...,4$, is a
para-K\"ahler structure on ${\mathcal G}(E)$.

For $K\in {\mathcal G}(E)$, let $K_1$ and $K_2$ be the
$g$-compatible paracomplex structures on $T$ defined by means of $K$
via (\ref{K1,2}). Then the map $K\to (K_1,K_2)$ is an isometry of
${\mathcal G}(E)$ onto $Z(T)\times Z(T)$. Moreover it sends a
tangent vector $V$ at $K\in {\mathcal G}(E)$ to the tangent vector
$(V_1,V_2)$, where
$$
V_1=(pr_{T}|E')\circ V \circ (pr_{T}|E')^{-1},\quad V_2=(pr_{T}|E'')\circ V \circ (pr_{T}|E'')^{-1}.
$$
Hence $K\to (K_1,K_2)$ is a bi-paraholomorphic map. Fixing an
orientation on the vector space $T$,  we see that ${\mathcal G}(E)$
has four connected components bi-paraholomorphically isometric to
$Z_{+}\times Z_{+}$, $Z_{+}\times Z_{-}$, $Z_{-}\times Z_{+}$,
$Z_{-}\times Z_{-}$.

\subsection{Generalized paracomplex structures and generalized metrics on manifolds}

A generalized almost paracomplex structure \cite {Wa, Va07} on a
smooth manifold $M$ is an endomorphism ${\mathcal K}$ of the bundle
$TM\oplus T^{\ast}M$ with ${\mathcal K}^2=Id$ which is
skew-symmetric with respect to the neutral metric
$$
<X+\alpha,Y+\beta>=\frac{1}{2}[\alpha(Y)+\beta(X)],\quad X,Y\in TM,\> \alpha,\beta\in T^{\ast}M.
$$
Clearly, this is an analog of the notion of a generalized almost
complex structure.

A generalized almost paracomplex structure is called integrable if
its Nijenhuis tensor vanishes, the latter being defined by means of
the Courant bracket instead of the Lie one. Recall that the Courant
bracket \cite{Cou} is defined by
$$
[X+\alpha,Y+\beta]=[X,Y]+{\mathcal L}_{X}\beta-{\mathcal L}_{Y}\alpha-\frac{1}{2}d(\imath_X\beta-\imath_Y\alpha),
$$
where $[X,Y]$ on the right hand-side is the Lie bracket, ${\mathcal
L}$ means the Lie derivative, and $\imath$ stands for the interior
product. Note that the Courant bracket is skew-symmetric like the
Lie bracket, but it does not satisfy the Jacobi identity.

An integrable generalized almost paracomplex structure is also
called generalized paracomplex structure.

The proof of \cite[Proposition 3.23]{Gu} contains the following important property of the Courant bracket.
\begin{prop}\label{Courant-B-transf}
If $\Theta$ is a $2$-form on $M$, then for every sections  $A=X+\alpha$ and $B=Y+\beta$ of $TM\oplus T^{\ast}M$
$$
[e^{\Theta}A,e^{\Theta}B]=e^{\Theta}[A,B]-\imath_{X}\imath_{Y}d\Theta.
$$
\end{prop}
Thus, if the form $\Theta$ is closed and ${\mathcal K}$ is a
generalized paracomplex structure,  the generalized paracomplex
structure $e^{\Theta}{\mathcal K}e^{-\Theta}$ is integrable exactly
when the structure ${\mathcal K}$ is integrable.

If ${\mathcal K}$ is a generalized paracomplex structure, let
${\mathcal L}^{\pm}$ be the subbundles corresponding to the
eigenvalue $\pm$ of ${\mathcal K}$. Then ${\mathcal K}$ is
integrable if and only if the space of sections of ${\mathcal
L}^{\pm}$ is closed under the Courant bracket \cite{Wa}.

\smallskip

\noindent {\bf Examples} \cite{Wa}.

\smallskip

\noindent {\bf 7}. The trivial  generalized almost paracomplex
structure ${\mathcal K}(X+\alpha)=X-\alpha$ is integrable.

\smallskip

\noindent {\bf 8}. The generalized almost  paracomplex structure
${\mathcal K}_{\omega}(X+\alpha)=\omega^{-1}(\alpha)+\omega(X)$
defined by means of a pre-symplectic form $\omega$ (non-degenerate
skew-symmetric $2$-form) is integrable if and only if $d\omega=0$,
i.e., $\omega$ is a symplectic form.

\smallskip

\noindent {\bf 9}. Let $\pi$ be a field of $2$-vectors on a
manifold. The generalized almost paracomplex structure ${\mathcal
K}_{\pi}(X+\alpha)=(X-\imath_{\alpha}\pi)-\alpha$ is integrable if
and only if $\pi$ is Poisson.

\smallskip

\noindent {\bf 10}. If $P$ is an almost paracomplex structure, the
generalized almost paracomplex structure ${\mathcal K}_P$ is
integrable if and only if $P$ is integrable.

\medskip

By definition, a generalized metric on a manifold $M$ is a subbundle
$E$ of $TM\oplus T^{\ast}M$ such that $rank\, E=dim\, M$, the
restriction of the metric
$<X+\alpha,Y+\beta>=\frac{1}{2}[\alpha(Y)+\beta(X)]$ to $E$ is of
signature $(r,s)$, and $E_x\cap T_x^{\ast}M=\{0\}$ for every $x\in
M$. Every such a bundle is uniquely determined by a metric $g$ of
signature $(r,s)$ and a skew-symmetric $2$-form $\Theta$ on $M$ such
that $E=\{X+g(X)+\Theta(X):~X\in TM\}$.

A generalized paracomplex structure ${\mathcal K}$ is called
compatible with $E$ if ${\mathcal K}E\subset E$. Clearly, in this
case, ${\mathcal K}E=E$ and ${\mathcal K}E^{\perp}=E^{\perp}$, the
orthogonal complement being taken with respect to the metric
$<.\,,\,.>$.

\smallskip

\noindent {\bf Examples}.

\smallskip

\noindent {\bf 11}. Let $\omega$ be a pre-symplectic form on a
manifold $M$, $dim\,M=2n$. Suppose there is an almost paracomplex
structure $L$ on $M$ compatible with $\omega$ in the sense that
$\omega(LX,LY)=-\omega(X,Y)$. Then $g(X,Y)=\omega(LX,Y)$ is a
neutral metric and the generalized paracomplex structure ${\mathcal
K}_{\omega}$ is compatible with the generalized metric
$\{X+g(X):~X\in TM\}$ (Example  4).

Note that an almost paracomplex structure compatible with $\omega$
exists if and only if the structure group of $M$ can be reduced from
$Sp(2n)$ to the paraunitari group $/U(n,{\mathbb A})$ \cite{EST, M}.
This is a restrictive condition and there are symplectic manifolds
which do not admit any compatible almost paracomplex structure, see
\cite[Sec. 2.5]{EST}, \cite[Sec. 2]{M}. Of course, any para-K\"ahler
manifold considered as a symplectic one possesses a compatible
paracomplex structure.

\smallskip

\noindent {\bf 12}. Other examples of generalized almost paracomplex
structures compatible with a generalized metric on a manifold can be
obtained by means of Examples 5 and 6 in an obvious way.

\subsection{Connections induced by  generalized
metrics on manifolds}

It has been shown by Hitchin \cite{Hit06,Hit10} that if
$E=\{X+g(X)+\Theta(X):~X\in T\}$ is a positive definite generalized
metric one can define a connection $\nabla^E$ on the bundle $E\to M$
preserving the metric $<.\,,.>|E$. Transferring the connection
$\nabla^E$ from the bundle $E$ to the bundle $TM$ via the
isomorphism $pr_{TM}|E:E\to TM$ we get a connection $\nabla$ on $TM$
preserving the metric $g$. Hitchin [ibid]  has proved that the
torsion $T$ of the connection $\nabla$ is totally skew-symmetric and
is given by
$$
g(T(X,Y),Z)=d\Theta(X,Y,Z),\quad X,Y,Z\in TM.
$$
Hitchin's proof also goes without any changes in the case of a
generalized metric of an arbitrary signature.

\medskip
\noindent {\bf Convention}. Henceforward, we consider only metrics
of neutral  signature $(n,n)$.

\medskip

\section{Generalized reflector spaces}

  Let $M$ be a smooth manifold of dimension $m=2n$ equipped with a generalized metric $E$ of signature $(n,n)$ determined
by a metric $g$ of signature $(n,n)$ and a $2$-form $\Theta$ on $M$. Denote by ${\mathcal G}={\mathcal G}(E)\to M$  the
bundle over $M$ whose fibre at a point $p\in M$ consists of all generalized paracomplex structures on $T_pM$ compatible
with the generalized metric $E_p$, the fibre of $E$ at $p$. We call ${\mathcal G}$  generalized reflector space of the
generalized Riemannian manifold $(M,E)$.

Set $E'=E$ and $E''=E^{\perp}$, the orthogonal complement of $E$ in
$TM\oplus T^{\ast}M$ with respect to the metric $<.\,,.>$. Denote by
${\mathcal Z}(E')$ the bundle over $M$ whose fibre at a point $p\in
M$ is the manifold $Z(E'_p)$ of  paracomplex structure on the vector
space $E'_p$ compatible with the neutral metric $g'=<.\,,.>|E'$.
Define a bundle ${\mathcal Z}(E'')$ in a similar way. Then
${\mathcal G}$ is diffeomorphic to  the product bundle ${\mathcal
Z}(E')\times {\mathcal Z}(E'')$ by the map ${\mathcal G}_p\ni K\to
(K|E',K|E'')$.

Suppose we are given metric connections $D'$ and $D''$ on $E'$ and $E''$, respectively,  and let $D=D'\oplus D''$ be the
connection on $E'\oplus E''=TM\oplus T^{\ast}M$ determined by $D'$ and $D''$.

 The bundle ${\mathcal Z}(E')$ is a subbundle of the vector bundle
$A(E')$ of $g'$-skew-symmetric endomorphisms of $E'$, and similarly
for ${\mathcal Z}(E'')$. Henceforth  we shall consider the bundle
${\mathcal G}\cong {\mathcal Z}(E')\times {\mathcal Z}(E'')$ as a
subbundle of the vector bundle $\pi:A(E')\oplus A(E'')\to M$. The
connection on $A(E')\oplus A(E'')$ induced by the connection
$D=D'\oplus D''$ on $E'\oplus E''$ will again be denoted by $D$. It
is easy to see that the horizontal space of $A(E')\oplus A(E'')$
with respect to $D$ at every point of ${\mathcal G}$ is tangent to
${\mathcal G}$ (cf. the next section). Thus the connection $D$ gives
rise to a splitting ${\mathcal V}\oplus {\mathcal H}$ of the tangent
bundle of the manifold ${\mathcal G}$ into vertical and horizontal
parts. Then, following the standard twistor construction,  we can
define four generalized almost paracomplex structures ${\mathscr
K}_{\varepsilon}$ on the manifold ${\mathcal G}$.

 The vertical space ${\mathcal V}_K$ of ${\cal G}$
at a point $K\in {\mathcal G}$ is the tangent space at $K$ of the
fibre through this point. This fibre is the manifold ${\mathcal
G}(E_{\pi(K)})$, which admits four paracomplex structures ${\mathscr
P}_{\varepsilon}$ defined in Sec.~\ref{CGPCS}. We define ${\mathscr
K}_{\varepsilon}|({\mathcal V}_K\oplus {\mathcal V}_K^{\ast})$ to be
the generalized paracomplex structure determined by the paracomplex
structure ${\mathscr P}_{\varepsilon}$. Thus
$$
{\mathscr K}_{\varepsilon}={\mathscr P}_{\varepsilon}~ \rm{on}~ {\mathcal V}_K,\quad   {\mathscr K}_{\varepsilon}=
-{\mathscr P}_{\varepsilon}^{\ast} ~\rm{on}~ {\mathcal V}_K^{\ast},\quad \varepsilon=1,2,3,4.
$$

 The horizontal space ${\cal H}_Ê$ is isomorphic  to the tangent space $T_{\pi(K)}M$ via the
differential $\pi_{\ast Ê}$. If $\pi_{{\mathcal H}}$ is the
restriction of $\pi_{\ast}$ to ${\cal H}$, the image of every $A\in
T_pM\oplus T_p^{\ast}M$ under the map $\pi_{{\mathcal H}}^{-1}\oplus
\pi_{{\mathcal H}}^{\ast}$ will be denoted by $A^h$. Thus, for $K\in
{\mathcal G}$, $X\in T_{\pi(K)}M$, and $\omega\in
T_{\pi(K)}^{\ast}M$, we have $\omega^h_K(X^h_Ê)=\omega_{\pi(K)}(X)$.
The elements of ${\mathcal H}_K^{\ast}$, resp. ${\mathcal
V}_K^{\ast}$, will be considered as $1$-forms on $T_K{\cal G}$
vanishing on ${\mathcal V}_K$, resp. ${\mathcal H}_K$. Then the
metrics $<.\,,.>$ on the spaces $T_K{\mathcal G}\oplus
T_K^{\ast}{\mathcal G}$, $T_{\pi(K)}M\oplus T_{\pi(K)}^{\ast}M$,
${\mathcal V}_K\oplus {\mathcal V}_K^{\ast}$ are related by
$$
<A^h+{\mathscr U},B^h+{\mathscr W}>=<A,B>+<{\mathscr U},{\mathscr W}>,
$$
$A,B\in T_{\pi(K)}M\oplus T_{\pi(K)}^{\ast}M$,~ ${\mathscr
U},{\mathscr W}\in {\mathcal V}_K\oplus {\mathcal V}_K^{\ast}$.

We define a generalized paracomplex structure ${\mathscr K}$ on the
vector space ${\mathcal H}_K\oplus {\mathcal H}_K^{\ast}$ as the
lift of the endomorphism $K$  of $T_{\pi(K)}M\oplus
T_{\pi(K)}^{\ast}M$ by the isomorphism $\pi_{{\mathcal H}}\oplus
(\pi_{{\mathcal H}}^{-1})^{\ast}:{\mathcal H}_K\oplus {\mathcal
H}_K^{\ast}\to T_{\pi(K)}M\oplus T_{\pi(K)}^{\ast}M$:
$$
{\mathscr K}A^h_K=(KA)^h_K,\quad A\in T_{\pi(K)}M\oplus T^{\ast}_{\pi(K)}M.
$$

Now, set ${\mathscr K}_{\varepsilon}={\mathscr K}$ on ${\mathcal
H}\oplus {\mathcal H}^{\ast}$.

\begin{prop}
The generalized paracomplex structure ${\mathcal K}_{\varepsilon}$
is not a $B$-transform of the trivial generalized almost paracomplex
structure ${\mathcal K}$ or the generalized almost paracomplex
structures ${\mathcal K}_{\omega}$, ${\mathcal K}_{\pi}$, ${\mathcal
K}_{P}$ yielded by a pre-symplectic form $\omega$, a field of
$2$-vectors $\pi$, a product structure $P$ on ${\mathcal G}$,
respectively.
\end{prop}

\begin{proof}

Note first that the restriction of any $B$-transform $e^B$ to the
cotangent bundle of a  manifold is the identity.

Let $E=\{X+g(X)+\Theta(X):~X\in TM\}$. Fix a point $p\in M$, and
take a paracomplex structure $L$ on $T_pM$ compatible with the
metric $g_p$; in particular, $L\neq \pm Id$. Applying
Proposition~\ref{K-K1,2} with $K_1=K_2=L$, we get a generalized
paracomplex structure $K$ on $T_pM$ compatible with the generalized
metric $E_p$ and such that
$K(X+\alpha)=LX+\Theta(LX)+L^{\ast}(\Theta(X))-L^{\ast}\alpha$ for
$X\in T_pM$ and $\alpha\in T_p^{\ast}M$.

\smallskip

Suppose ${\mathcal K}=e^{-B}{\mathcal K}_{\varepsilon}e^B$ for some
skew-symmetric $2$-form $B$ on ${\mathcal G}$. Then, for $\alpha\in
T_p^{\ast}M$, denoting the $T_pM$ and $T_p^{\ast}M$-components of
$K\alpha$ by $(K\alpha)_{T}$ and $(K\alpha)_{T^{\ast}}$,
$$
-\alpha^h_K={\mathcal
K}(\alpha^h_K)=e^{-B}((K\alpha)^h_K)=((K\alpha)_T)^h_K+((K\alpha)_{T^{\ast}})^h_K-\imath_{((K\alpha)_T)^h_K}B
=-(L^{\ast}\alpha)^h_K.
$$
Hence $L^{\ast}\alpha=\alpha$ for every $\alpha\in T_p^{\ast}M$.
Then $\alpha(LX)=\alpha(X)$ for $X\in T_pM$, which implies $LX=X$.
Thus, $L=Id$, a contradiction

\smallskip

Let ${\mathcal K}_{\omega}=e^{-B}{\mathcal K}_{\varepsilon}e^B$, and
take   $\Upsilon\in{\mathcal V}_K^{\ast}$, $\Upsilon\neq 0$. Then
$T{\mathcal G}\ni\omega^{-1}(\Upsilon)=e^{-B}{\mathcal
K}_{\varepsilon}(\Upsilon)=-e^{-B}{\mathscr
P}_{\varepsilon}^{\ast}(\Upsilon)=-{\mathscr
P}_{\varepsilon}^{\ast}(\Upsilon)\in T^{\ast}{\mathcal G}$. It
follows ${\mathscr P}_{\varepsilon}^{\ast}(\Upsilon)=0$. Hence,
$\Upsilon=0$, a contradiction.

\smallskip

Next, suppose ${\mathcal K}_{\pi}=e^{-B}{\mathcal
K}_{\varepsilon}e^B$. Then,
$$
-\imath_{\alpha^h_K}\pi-\alpha^h_K=-(L^{\ast}\alpha)^h_K.
$$
Hence $L^{\ast}\alpha=\alpha$, a contradiction.

\smallskip

Finally, suppose ${\mathcal K}_{P}=e^{-B}{\mathcal
K}_{\varepsilon}e^B$. Set $\omega(X,Y)=g(X,LY)$, and apply
Proposition~\ref{K-K1,2} with $K_1=-L$, $K_2=L$; so, under the
notation of this proposition, $-\omega_1=\omega_2=\omega$. In this
way, we obtain a generalized paracomplex structure $Q\in {\mathcal
G}(E_p)$ such that
$Q(X+\alpha)=-(\omega^{-1}\circ\theta)(X)+\omega^{-1}(\alpha)-(\Theta\circ\omega^{-1}\circ\Theta)(X)
+(\Theta\circ\omega^{-1})(\alpha)+\omega(X)$. Then
$$
-P^{\ast}\alpha^h_Q=(\omega^{-1}(\alpha))^h_Q
+((\Theta\circ\omega^{-1})(\alpha))^h_Q-\imath_{(\omega^{-1}(\alpha))^h_Q}B.
$$
This implies $\omega^{-1}(\alpha)=0$, hence $\alpha=0$ fo every
$\alpha$.

\end{proof}

\section{Technicalities}

Let $(\mathscr{U},x_1,...,x_{2n})$ be a local coordinate system of
$M$ and $\{Q_1',...,Q_{2n}'\}$, $\{Q_1'',...,Q_{2n}''\}$ orthonormal
frames of $E'$ and $E''$ on $\mathscr{U}$ with respect to the
metrics $g'=<.\,,.>|E'$ and $g''=-<.\,,.>|E''$, respectively, such
that $||Q'_{\alpha}||^2_{g'}=-||Q'_{n+\alpha}||^2_{g'}=1$,
$||Q''_{\alpha}||^2_{g''}=-||Q''_{n+\alpha}||^2_{g''}=1$,
$\alpha=1,\dots, n$. Define sections $S_{ij}'$, $S_{ij}''$, $1\leq
i,j\leq {2n}$, of $A(E')$ and $A(E'')$ by the formulas
\begin{equation}\label{eq Sij}
S_{ij}'Q_k'=\delta_{ik}||Q_j'||^2_{g'}Q_j' - \delta_{kj}||Q_i||^2_{g'}Q_i',\quad
S_{ij}''Q_k''=\delta_{ik}||Q_j''||^2_{g''}Q_j'' - \delta_{kj}||Q_i||^2_{g''}Q_i''.
\end{equation}
Then $S_{ij}'$ and  $S_{ij}''$ with $i<j$ form orthonormal frames of $A(E')$ and $A(E'')$ with respect to the metrics
$G'$ and $G''$ defined by
$$
G'(a',b')=\displaystyle{-\frac{1}{2}}Trace\,(a'\circ b').\quad
G''(a'',b'')=\displaystyle{-\frac{1}{2}}Trace\,(a''\circ b'')
$$
for $a',b'\in A(E')$ and $a'',b''\in A(E'')$. Thus, every $a'\in
A(E')$ can be written as
$$
a'=\sum\limits_{i<j}G'(a',S_{ij}')||S_{ij}'||^2_{G'}S_{ij}'
$$
where $||S_{ij}'||^2_{G'}=\frac{1}{2}(||Q_i'||^2_{g'}+||Q_j'||^2_{g'})$. Similarly for $a''\in A(E'')$.

\smallskip

 For $a=(a',a'')\in A(E')\oplus A(E'')$, set
\begin{equation}\label{coord}
\begin{array}{c}
\tilde x_{i}(a)=x_{i}\circ\pi(a),\\[6pt]
y_{kl}'(a)=G'(a',S_{kl}'\circ\pi(a))||S_{kl}'\circ\pi(a)||^2_{G'},\quad
y_{kl}''(a)=G''(a'',S_{kl}''\circ\pi(a))||G_{kl}''\circ\pi(a)||^2_{G''}.
\end{array}
\end{equation}
Then $(\tilde x_{i},y_{jk}', y_{jk}'')$, $1\leq i\leq 2n$, $1\leq j < k\leq 2n$,  is a local coordinate system on the
total space of the bundle $A(E')\oplus A(E'')$. Note also that
$$
y_{kl}'(a)=g'(a'(Q_k'(p)),Q_l'(p)),\quad
y_{kl}''(a)=g''(a''(Q_k''(p)),Q_l''(p)),\quad p=\pi(a),
$$
thus $a'(Q_k'(p))=\sum\limits_{l=1}^{2n}
y_{kl}(a)||Q_l'(p)||^2_{g'}Q_l'(p)$ and similar for $a''(Q_k''(p))$,

Let
\begin{equation}\label{V}
V=\sum_{j<k}[v_{jk}'\frac{\partial}{\partial
y_{jk}'}(K)+v_{jk}''\frac{\partial}{\partial y_{jk}''}(K)]
\end{equation}
be a vertical vector of ${\cal G}$ at a point $K$.  It is convenient to set $v_{ij}'=-v_{ji}'$, $v_{ij}''=-v_{ji}''$  and
$y_{ij}'=-y_{ji}'$, $y_{ij}''=-y_{ji}''$ for $i\geq j$, $1\leq i,j\leq {2n}$. Then the endomorphism $V$ of $T_{p}M\oplus
T_{p}^{\ast}M$, $p=\pi(K)$, is determined by
$$
VQ_i'=\sum_{j=1}^{2n} v_{ij}'||Q_j'||^2_{g'}Q_j',\quad VQ_i''=\sum_{j=1}^{2n} v_{ij}''||Q_j''||^2_{g''}Q_j''.
$$
Hence 
\begin{equation}\label{cal J/ver}
{\mathscr K}_{\varepsilon}V=(-1)^{\varepsilon+1}\sum_{j<k} \sum_{s=1}^{2n}[\pm
v_{js}'y_{sk}'||Q_s'||^2\frac{\partial}{\partial y_{jk}'}+ v_{js}''y_{sk}''||Q_s''||^2\frac{\partial}{\partial y_{jk}''}]
\end{equation}
where the plus sign corresponds to $\varepsilon=1,4$ and the minus sign to $\varepsilon=2,3$.

\smallskip

 Note also that, for every $A\in T_{p}M\oplus T_{p}^{\ast}M$,

\begin{equation}\label{cal K/hor}
\begin{array}{l}
 A^h=\sum\limits_{i=1}^{2n}[(\big(<A,Q_i'>||Q_i'||^2\big)\circ\pi)Q_i^{'\,h}-(\big(<A,Q_i''>||Q_i''||^2\big)\circ\pi)Q_i^{''\,h}]\\[8pt]
{\mathscr K}A^h=\sum\limits_{i,j=1}^{2n}[(\big(<A,Q_i'>||Q_i'||^2||Q_j'||^2\big)\circ\pi)y_{ij}'Q_j^{'\,h}\\[6pt]
\hspace{4.8cm}-(\big(<A,Q_i''>||Q_i''||^2||Q_j''||^2\big)\circ\pi)y_{ij}''Q_j^{''\,h}].
\end{array}
\end{equation}

For each vector field
$$X=\sum_{i=1}^{2n} X^{i}\frac{\partial}{\partial x_i}$$
on $\mathscr{U}$, the horizontal lift $X^h$ on $\pi^{-1}(\mathscr{U})$ is given by
\begin{equation}\label{Xh}
\begin{array}{c}
X^{h}=\displaystyle{\sum_{l=1}^{2n} (X^{l}\circ\pi)\frac{\partial}{\partial\tilde
x_l}}\\[8pt]
- \displaystyle{\sum_{i<j}\sum_{k<l}
y_{kl}'[\big(G'(D_{X}S_{kl}',S_{ij}')||S_{ij}||^2_{G'}\big)\circ\pi]\frac{\partial}{\partial y_{ij}'}}\\[8pt]
-\displaystyle{\sum_{i<j}\sum_{k<l}
y_{kl}''[\big(G''(D_{X}S_{kl}'',S_{ij}'')||S_{ij}''|^2_{G''}\big)\circ\pi]\frac{\partial}{\partial y_{ij}''}}.
\end{array}
\end{equation}

Let $a=(a',a'')\in A(E')\oplus A(E'')$. Denote  the fiber of $A(E')$
at the point $\pi(a)$ by $A(E'_{\pi(a)})$, and similarly for
$A(E''_{\pi(a)})$. Then (\ref{Xh}) implies that, under the standard
identification of $T_{a}(A(E'_{\pi(a)})\oplus A(E''_{\pi(a)}))$ with
the vector space $A(E'_{\pi(a)})\oplus A(E''_{\pi(a)})$,
\begin{equation}\label{[XhYh]}
[X^{h},Y^{h}]_{a}=[X,Y]^h_a + R(X,Y)a
\end{equation}
where $R(X,Y)a=R(X,Y)a'+R(X,Y)a''$ is the curvature of the
connection $D$ (for the curvature tensor {\it we adopt the following
definition}: $R(X,Y)=D_{[X,Y]}- [D_{X},D_{Y}]$). Note also that
(\ref{V}) and (\ref{Xh}) imply the well-known fact that
\begin{equation}\label{[V,Xh]}
[V,X^h]~\rm{is~ a~ vertical~ vector~ field}.
\end{equation}

\smallskip

\noindent {\it Notation}. Let $K\in {\mathcal G}$ and $p=\pi(K)$.
Take  orthonormal bases $\{a_1',...,a_{2n}'\}$,
$\{a_1'',...,a_{2n}''\}$ of $E_p'$ and $E_p''$ such that
$a_{2l}'=Ka_{2l-1}'$, $a_{2l}''=Ka_{2l-1}''$ for $l=1,...,n$. Let
$\{Q_i'\}$, $\{Q_i''\}$, $i=1,...,2n$, be  orthonormal frames of
$E'$, $E''$ near the point $p$ such that
$$
Q_i'(p)=a_i',\> Q_i''(p)=a_i'' ~\mbox { and }~ D\, Q_i'|_p=0,\> D\,Q_i''|_p=0, \quad i=1,...,2n.
$$
Define a section $S=(S',S'')$ of $A(E')\oplus A(E'')$  setting
$$
S'Q_{2l-1}'=Q_{2l}',\quad S''Q_{2l-1}''=Q_{2l}'' ,\quad S'Q_{2l}'=Q_{2l-1}',\quad S''Q_{2l}''=Q_{2l-1}''
$$
$l=1,...,n$. Then,
$$
S(p)=K, ~ D S|_p=0.
$$
In particular $X^h_K=S_{\ast}X$ for every $X\in T_pM$.

Clearly, the section $S$ takes its values in ${\mathcal G}$, hence
the horizontal space of $A(E')\oplus A(E'')$ with respect to the
connection $D$ at any $K\in {\mathcal G}$ is tangent to ${\cal G}$.

 Further on, given a smooth manifold $M$, the natural projections of $TM\oplus T^{\ast}M$
onto $TM$ and $T^{\ast}M$ will be denoted by $\pi_1$ and $\pi_2$,
respectively. The natural projections of ${\mathcal H}\oplus
{\mathcal H}^{\ast}$ onto ${\mathcal H}$ and ${\mathcal H}^{\ast}$
will also be denoted by $\pi_1$ and $\pi_2$ when this will not cause
confusion. Thus if $\pi_1(A)=X$ for $A\in TM\oplus T^{\ast}M$, then
$\pi_1(A^h)=X^h$ and similarly for $\pi_2(A)$ and $\pi_2(A^h)$.

  We shall use the above notations throughout the next sections.

\smallskip

Note that, although $DS|_p=0$, $D\pi_1(S)$ and $D\pi_2(S)$ may not vanish at the point $p$ since the connection $D$ may
not preserve $TM$ or $T^{\ast}M$.

\smallskip

\begin{lemma}\label {brackets}
If $A$ and $B$ are sections of the bundle $TM\oplus T^{\ast}M$ near $p$, then
\begin{enumerate}
\item[$(i)$] $[\pi_1(A^h),\pi_1(B^h)]_K=[\pi_1(A),\pi_1(B)]^h_K+R(\pi_1(A),\pi_1(B))K.$
\vspace{0.2cm}
\item[$(ii)$]  $[\pi_1(A^h),\pi_1({\mathscr
K}B^h)]_K=[\pi_1(A),\pi_1(SB)]^h_K+R(\pi_1(A),\pi_1(KB))K.$ \vspace{0.2cm}
\item[$(iii)$]
$[\pi_1({\mathscr K}A^h),\pi_1({\mathscr K}B^h)]_K=[\pi_1(SA),\pi_1(SB)]^h_K+R(\pi_1(KA),\pi_1(KB))K.$
\end{enumerate}
\end{lemma}

\begin{proof}
Identity $(i)$ is a paraphrase of (\ref{[XhYh]}).

To prove $(ii)$, set $X=\pi_1(A)$. By (\ref{Xh}), we have
$X_K^h=\displaystyle{\sum\limits_{l=1}^{2n}
X^l(p)\frac{\partial}{\partial\tilde x_l}(K)}$ since $D S_{kl}'|_p=D
S_{kl}''|_p=0$, $k,l=1,...,2m$. Then, using (\ref{cal K/hor}), we
get
\begin{equation}\label{XhBh}
\begin{array}{c}
[X^h,\pi_1({\mathscr K}B^h)]_K=\\[8pt]

\sum\limits_{i,j=1}^{2n}\big[<B,Q_i'>_p||Q_i'||^2_p||Q_j'||^2_p y_{ij}'(K)[X^h,\pi_1(Q_j')^h_K\\[6pt]

\hspace{4.5cm}+X_p(<B,Q_i'>||Q_i'||^2||Q_j'||^2)y_{ij}'(K)(\pi_1(Q_j'))^h_K\big]\\[8pt]

-\sum\limits_{i,j=1}^{2n}\big[<B,Q_i''>_p||Q_i''||^2_p||Q_j''||^2_py_{ij}''(K)[X^h,\pi_1(Q_j')^h]_K\\[6pt]

\hspace{4.5cm}- X_p(<B,Q_i''>||Q_i||^2||Q_j||^2)y_{ij}''(K)(\pi_1(Q_j''))^h_K\big].
\end{array}
\end{equation}
We also have
\begin{equation}\label{sB}
SB=\sum_{i,j=1}^{2n}\big[<B,Q_i'>||Q_i'||^2||Q_j'||^2(y_{ij}'\circ
S)Q_j'-<B,Q_i''>||Q_i''||^2||Q_j''||^2(y_{ij}''\circ S)Q_j''\big].
\end{equation}
Therefore
\begin{equation}\label{XB}
\begin{array}{c}
[X,\pi_1(SB)]_p=\\[8pt]

\sum\limits_{i,j=1}^{2n}\big[<B,Q_i'>_p||Q_i'||^2_p||Q_j'||^2_py_{ij}'(K)[X,\pi_1(Q_j')]_p\\[6pt]

\hspace{4.5cm}+X_p(<B,Q_i'>||Q_i'||^2||Q_j'||^2)y_{ij}'(K)(\pi_1(Q_j'))_p\big]\\[8pt]

-\sum\limits_{i,j=1}^{2n}\big[<B,Q_i''>||Q_i''||^2_p||Q_j''||^2_py_{ij}''(K)[X,\pi_1(Q_j'')]_p\\[6pt]

\hspace{4.5cm}+X_p(<B,Q_i''>||Q_i''||^2||Q_j''||^2)y_{ij}''(K)(\pi_1(Q_j''))_p\big].
\end{array}
\end{equation}

Now, formula $(i)$ follows from (\ref{XhBh}), (\ref{[XhYh]}) and
(\ref{XB}).

A similar computation gives $(iii)$.
\end{proof}

\smallskip

 For any (local) section $a=(a',a'')$ of $A(E')\oplus A(E'')$,
denote by $\widetilde a$ the vertical vector field on ${\mathcal G}$ defined by
\begin{equation}\label{eq tilde a}
\widetilde a_K=(a'_{\pi(K)}-(K|E'_{\pi(K)})\circ a'_{\pi(K)}\circ
(K|E'_{\pi(K)}, a''_{\pi(K)}-(K|E''_{\pi(K)}\circ a''_{\pi(K)}\circ
(K|E''_{\pi(K)})) .
\end{equation}
Let us note that, for every $K\in {\mathcal G}$, we can find
sections $a_1,...,a_s$, $s=2(n^2-n)$, of $A(E')\oplus A(E'')$ near
the point $p=\pi(K)$ such that $\widetilde a_1,...,\widetilde a_s$
form a basis of the vertical vector space at each point in a
neighbourhood of $K$.

\vspace{0.1cm}

\begin{lemma}\label {H-V brackets}
Let $K\in {\cal G}$ and let $a$ be a section of $A(E')\oplus A(E'')$ near the point $p=\pi(K)$. Then, for any section $A$
of the bundle $TM\oplus T^{\ast}M$ near $p$, we have (for the Lie brackets)
\begin{enumerate}
\item[$(i)$] $[\pi_1(A^h),\widetilde a]_K=(\widetilde{D_{\pi_1(A)}a})_K$
\vspace{0.2cm}
\item[$(ii)$] $[\pi_1(A^h),{\mathscr K}_{\varepsilon}\widetilde a]_K=
{\mathscr P}_{\varepsilon}(\widetilde{D_{\pi_1(A)}a})_K$ \vspace{0.2cm}
\item[$(iii)$] $[\pi_1({\mathscr K}A^h),\widetilde
a]_K=(\widetilde{D_{\pi_1(KA)}a})_K-(\pi_1(\widetilde a(A)))^h_K$ \vspace{0.2cm}
\item[$(iv)$] $[\pi_1({\mathscr K}A^h),{\mathscr K}_{\varepsilon}\widetilde
a]_K={\mathscr P}_{\varepsilon}
(\widetilde{D_{\pi_1(KA)}a})_K-(\pi_1(({\mathscr
P}_{\varepsilon}\widetilde a)(A)))^h_K$
\end{enumerate}
\end{lemma}

\begin{proof}
Let
$$
a'Q_i'=\sum_{j=1}^{2n}a_{ij}'||Q_i'||^2Q_j',\quad
a''Q_i''=\sum_{j=1}^{2n}a_{ij}''||Q_i''||^2Q_j'',\quad i=1,...,2n.
$$
Then, in the local coordinates of $A(E')\oplus A(E'')$ introduced above,
$$
\widetilde a=\sum_{i<j}[\widetilde a_{ij}'\frac{\partial}{\partial y_{ij}'}+\widetilde a_{ij}''\frac{\partial}{\partial
y_{ij}''}]
$$
where
\begin{equation}\label{tilde aij}
\begin{array}{c}
\widetilde
a_{ij}'=a_{ij}'\circ\pi+\sum\limits_{k,l=1}^{2n}y_{ik}'(a_{kl}'\circ\pi)y_{lj}'(\big(||Q_k'||^2||Q_l'||^2\big)\circ\pi),\\[8pt]
a_{ij}''=a_{ij}''\circ\pi+\sum\limits_{k,l=1}^{2n}y_{ik}''(a_{kl}''\circ\pi)y_{lj}''(\big(||Q_k''||^2||Q_l''||^2\big)\circ\pi).
\end{array}
\end{equation}
Let us also note that for every vector field $X$ on $M$ near the
point $p$,  by (\ref{Xh}),
$$
\begin{array}{c}
X_K^h=\displaystyle{\sum_{i=1}^{2n} X^i(p)\frac{\partial}{\partial\tilde
x_i}(K)},\\[8pt]
\displaystyle{[X^h,\frac{\partial}{\partial y_{ij}'}]_K=[X^h,\frac{\partial}{\partial y_{ij}''}]_K}=0
\end{array}
$$
since $D S_{ij}'|_p=D S_{ij}''|_p=0$. Moreover,
\begin{equation}\label{Xa_ij}
(D_{X_p}a')Q_i'=\sum_{j=1}^{2n}X_p(a_{ij}')||Q_j'||^2Q_j', \quad
(D_{X_p}a'')Q_i''=\sum_{j=1}^{2n}X_p(a_{ij}'')||Q_j''||^2Q_j''
\end{equation}
since $D Q_i'|_p=D Q_i''|_p=0$. Now the lemma follows by simple computations making use of (\ref{cal J/ver}) and
(\ref{cal K/hor}).
\end{proof}

\begin{lemma}\label {Lie deriv}
Let $A$ and $B$ be sections of the bundle $TM\oplus T^{\ast}M$ near $p$, and let $Z\in T_pM$, $W\in {\cal V}_K$. Then
\begin{enumerate}
\item[$(i)$]
\hspace{0.2cm}$({\mathcal L}_{\pi_1(A^h)}{\pi_2(B^h)})_K=({\mathcal L}_{\pi_1(A)}{\pi_2(B)})^h_K.$

\vspace{0.3cm}

\item[$(ii)$]
\hspace{0.2cm}$({\mathcal L}_{\pi_1(A^h)}{\pi_2({\mathscr K}B^h)})_K=({\mathcal L}_{\pi_1(A)}{\pi_2(SB)})^h_K.$

\vspace{0.3cm}

\item[$(iii)$]
$
\begin{array}{lll}
({\mathcal L}_{\pi_1({\mathscr K}A^h)}\pi_2(B^h))_K(Z^h+W)=\\[6pt]
({\mathcal L}_{\pi_1(SA)}\pi_2(B))^h_K(Z^h)+(\pi_2(B))_p(\pi_1(WA)).
\end{array}
$

\vspace{0.3cm}

\item[$(iv)$]
$
\begin{array}{lll}
({\mathcal L}_{\pi_1({\mathscr K}A^h)}\pi_2({\mathscr
K}B^h))_K(Z^h+W)=\\[6pt]
({\mathcal L}_{\pi_1(SA)}\pi_2(SB))^h_K(Z^h)+(\pi_2(KB))_p(\pi_1(WA)).
\end{array}
$
\end{enumerate}
\end{lemma}

\begin{proof}
Formula $(i)$  follows from (\ref{[XhYh]}) and (\ref{[V,Xh]});
$(ii)$ is a consequence of $(i)$, (\ref{cal K/hor}) and (\ref{sB}).
A simple computations involving (\ref{cal K/hor}), (\ref{[XhYh]}),
(\ref{[V,Xh]}) and (\ref{sB}) gives formula $(iii)$; $(iv)$ follows
from $(iii)$, (\ref{cal K/hor}) and (\ref{sB}).
\end{proof}

\vspace{0.1cm}


\begin{lemma}\label {Half-diff}
Let $A$ and $B$ be sections of the bundle $TM\oplus T^{\ast}M$ near
$p$. Let $Z\in T_pM$ and $W\in {\cal V}_K$. Then
\begin{enumerate}
\item[$(i)$]
\hspace{0.2cm}$(d\>\imath_{\pi_1(A^h)}\pi_2(B^h))_K=(d\>\imath_{\pi_1(A)}\pi_2(B))^h_K$

\vspace{0.3cm}

\item[$(ii)$]
$
\begin{array}{lll}
(d\>\imath_{\pi_1(A^h)}\pi_2({\mathscr K} B^h))_K(Z^h+W)= \\[6pt]
(d\>\imath_{\pi_1(A)}\pi_2(SB))^h_K(Z^h)+ (\pi_2(WB))_p(\pi_1(A))
\end{array}
$

\vspace{0.3cm}

\item[$(iii)$]
$
\begin{array}{lll}
(d\>\imath_{\pi_1({\mathscr K}A^h)}\pi_2(B^h))_K(Z^h+W)= \\[6pt]
(d\>\imath_{\pi_1(SA)}\pi_2(B))^h_K(Z^h)+ (\pi_2(B))_p(\pi_1(WA))
\end{array}
$

\vspace{0.3cm}

\item[$(iv)$]
$
\begin{array}{lll}
(d\>\imath_{\pi_1({\mathscr K}A^h)}\pi_2({\mathscr K}B^h))_K(Z^h+W)=\\[6pt]
(d\>\imath_{\pi_1(SA)}\pi_2(SB))^h_K(Z^h)+ (\pi_2(WB))_p(\pi_1(KA))+ (\pi_2(KB))_p(\pi_1(WA)
\end{array}
$

\end{enumerate}
\end{lemma}

\begin{lemma}\label {H-V Lie der&diff}
Let $A$ be a section of the bundle $TM\oplus T^{\ast}M$ and $V$ a vertical vector field on ${\cal G}$. Then
\begin{enumerate}
\item[$(i)$] ${\mathcal L}_{V}\pi_2(A^h)=0$; \hspace{0.2cm}$\imath_{V}\pi_2(A^h)=0$.
\vspace{0.2cm}
\item[$(ii)$] ${\mathcal L}_{V}\pi_2({\mathscr K}A^h)=\pi_2((VA)^h)$;
\hspace{0.2cm} $\imath_{V}\pi_2({\mathscr K}A^h)=0$.
\end{enumerate}
\end{lemma}

\vspace{0.1cm}

\noindent {\it Notation}. Let $K\in {\mathcal G}$. For any fixed $\varepsilon=1,...,4$, take a basis
$\{U_{2t-1}^{\varepsilon},U_{2t}^{\varepsilon}={\mathscr K}_{\varepsilon}U_{2t-1}^{\varepsilon}\}$, $t=1,...,n^2-n$, of
the vertical space ${\mathcal V}_K$.  Let $a_{2t-1}^{\varepsilon}$ be sections of $A(E')\oplus A(E'')$ near the point
$p=\pi(K)$ such that $a_{2t-1}^{\varepsilon}(p)=U_{2t-1}^{\varepsilon}$ and $D a_{2t-1}^{\varepsilon}|_p=0$. Define
vertical vector fields $\widetilde a_{2t-1}^{\varepsilon}$ by (\ref{eq tilde a}). Then $\{\widetilde
a_{2t-1}^{\varepsilon},{\cal J}_{\varepsilon}\widetilde a_{2t-1}^{\varepsilon}\}$, $t=1,...,n^2-n$, is a frame of the
vertical bundle on ${\mathcal G}$ near the point $K$. Denote by $\{\beta_{2t-1}^{\varepsilon},\beta_{2t}^{\varepsilon}\}$
the dual frame of the bundle ${\mathcal V}^{\ast}$. Then $\beta_{2t}^{\varepsilon}={\mathscr
K}_{\varepsilon}\beta_{2t-1}^{\varepsilon}$.

\vspace{0.1cm}

  Under these notations, we have the following.

\begin{lemma}\label {H-Vstar Lie der}
Let $A$ be a section of the bundle $TM\oplus T^{\ast}M$ near the
point $p=\pi(K)$. Then for every $Z\in T_pM$, $s,r=1,....,2(n^2-n)$
and $\varepsilon=1,...,4$, we have
\begin{enumerate}
\item[$(i)$]$({\mathcal L}_{\pi_1(A^h)}\beta_s^{\varepsilon})_K(Z^h+U_r^{\varepsilon})=
-\beta_s^{\varepsilon}(R(\pi_1(A),Z)K).$ \vspace{0.2cm}
\item[$(ii)$]$({\mathcal L}_{\pi_1({\mathscr K}A^h)}\beta_s^{\varepsilon})_J(Z^h+U_r^{\varepsilon})
=-\beta_s^{\varepsilon}(R(\pi_1(KA),Z)K).$
\vspace{0.2cm}
\item[$(iii)$]$({\mathcal L}_{\pi_1(A^h)}{\mathscr K}_{\varepsilon}\beta_s^{\varepsilon})_J(Z^h+U_r^{\varepsilon})=-({\mathscr K}_{\varepsilon}\beta_s^{\varepsilon})(R(\pi_1(A),Z)K).$
\vspace{0.2cm}
\item[$(iv)$]$({\mathcal L}_{\pi_1({\mathscr
K}A^h)}{\mathscr K}_{\varepsilon}\beta_s^{\varepsilon})_J(Z^h+U_r^{\varepsilon})=-({\mathscr
K}_{\varepsilon}\beta_s^{\varepsilon})(R(\pi_1(KA),Z)K).$
\end{enumerate}
\end{lemma}

\begin{proof}
By (\ref{[XhYh]}), if $X=\pi_1(A)$,
$$
({\mathcal
L}_{\pi_1(A^h)}\beta_s^{\varepsilon})_K(Z^h+U_r^{\varepsilon})=
-\beta_s^{\varepsilon}(R(X,Z)K)-\frac{1}{2}\beta_s^{\varepsilon}([X^h,\widetilde
a_{r}^{\varepsilon}]_K).
$$
By  Lemma~\ref{H-V brackets},
$$
\begin{array}{c}
[X^h,a_{2t-1}^{\varepsilon}]_K=(\widetilde{D_{X}a_{2t-1}^{\varepsilon}})_K=0,\\[6pt]
[X^h,a_{2t}^{\varepsilon}]_K=[X^h,{\mathscr
K}_{\varepsilon}\widetilde
a_{2t-1}^{\varepsilon}]_K=K_{\varepsilon}(\widetilde{D_{X}a_{2t-1}^{\varepsilon}})_K=0
\end{array}
$$
since $D a_{2t-1}^{\varepsilon}|_p=0$. This proves the first
identity of the lemma. To prove the second one, we note that if $f$
is a smooth function on ${\mathcal G}$ and $Y$ is a vector field on
$M$, $ ({\mathcal
L}_{fY^h}\beta_s^{\varepsilon})_K(Z^h+U_r^{\varepsilon})=f({\mathcal
L}_{Y^h}\beta_s^{\varepsilon})_K(Z^h+U_r^{\varepsilon}) $ since
$\beta_s^{\varepsilon}(Z^h)=0$. Now, $(ii)$ follows from (\ref{cal
K/hor}) and the first identity of the lemma. Identities $(iii)$ and
$(iv)$ are straightforward consequences of $(i)$ and $(ii)$,
respectively, since ${\mathscr
K}_{\varepsilon}\beta_{2t-1}^{\varepsilon}=\beta_{2t}^{\varepsilon}$,
${\mathscr
K}_{\varepsilon}\beta_{2t}^{\varepsilon}=\beta_{2t-1}^{\varepsilon}$,
$t=1,...,n^2-n$.
\end{proof}


\section{The Nijenhuis tensor}

\noindent {\it Notation}. We denote the Nijenhuis tensor of
${\mathscr K}_{\varepsilon}$ by ${\mathscr N}_{\varepsilon}$,
$\varepsilon=1,2,3,4$,:
$$
{\mathscr N}_{\varepsilon}({\mathcal E},{\mathcal F})=[{\mathcal E},{\mathcal F}]+[{\mathscr K}_{\varepsilon}{\mathcal
E},{\mathscr K}_{\varepsilon}{\mathcal F}]-{\mathscr K}_{\varepsilon}[{\mathscr K}_{\varepsilon}{\mathcal E},{\mathcal
F}]-{\mathscr K}_{\varepsilon}[{\mathcal E},{\mathscr K}_{\varepsilon}{\mathcal F}],\quad {\mathcal E},{\mathcal F}\in
T{\mathcal G}.
$$
Moreover, given $K\in{\mathcal G}$ and $A,B\in T_{p}M\oplus T_{p}^{\ast}M$, $p=\pi(K)$, we define $1$-forms on ${\mathcal
V}_K$  setting
$$
\omega^{\varepsilon}_{A,B}(W)=<({\mathscr P}_1W-{\mathscr
P}_{\varepsilon}W)(A),B>-<({\mathscr P}_1W-{\mathscr
P}_{\varepsilon}W)(B),A>, \quad W\in{\mathcal V}_J.
$$

Also, let $S$ be a section of ${\cal G}$ in a neighbourhood of the point $p=\pi(K)$  such that $S(p)=K$ and $D S|_p=0$
($S$ being considered as a section of $A(E')\oplus A(E'')$).

Denote the projection operators onto the horizontal and vertical components by ${\cal H}\oplus {\cal H}^{\ast}$ and
${\mathcal V}\oplus {\mathcal V}^{\ast}$,

\vspace{0.1cm}
\begin{prop}\label {Nijenhuis}
Let $K\in {\cal G}$, $A,B\in T_{\pi(K)}M\oplus T_{\pi(K)}^{\ast}M$, $V,W\in {\mathcal V}_K$, $\varphi,\psi\in {\mathcal
V}_K^{\ast}$. Then:

\begin{enumerate}

\item[$(i)$] Extending $A$ and $B$ to sections of $TM\oplus T^{\ast}M$, we have\ \\

$ ({\cal H}\oplus {\cal H}^{\ast}){\mathscr N}_{\varepsilon}(A^h,B^h)_K=([A,B] +[SA,SB]-S[A,SB]-S[SA,B])^h_K. $

\vspace{0.1cm}

\item[$(ii)$] \ \\
$$
\begin{array}{c}
({\mathcal V}\oplus {\mathcal V}^{\ast}){\mathscr N}_{\varepsilon}(A^h,B^h)_K
= R(\pi_1(A),\pi_1(B))K +R(\pi_1(KA),\pi_1(KB))K \\[8pt]

\hspace{4.5cm}-{\mathscr P}_{\varepsilon}R(\pi_1(KA),\pi_1(B))K - {\mathscr P}_{\varepsilon}R(\pi_1(A),\pi_1(KB))K \\[8pt]

-\omega^{\varepsilon}_{A,B},
\end{array}
$$

\vspace{0.1cm}

\item[$(iii)$] \ \\
${\mathscr N}_{\varepsilon}(A^h,V)_K=(-({\mathscr P}_{\varepsilon}V)A+({\mathscr P}_1V)A)^{h}_K.$

\vspace{0.1cm}

\item[$(iv)$] \ \\
${\mathscr N}_{\varepsilon}(A^h,\varphi)_K\in {\cal H}_K\oplus {\cal H}_K^{\ast}$ ~ and \ \\
              \ \\
$ <\pi_{\ast}N_{\varepsilon}(A^h,\varphi)_K,B>=-\displaystyle{\frac{1}{2}}\varphi({\mathcal
V}N_{\varepsilon}(A^h,B^h)_K).$

\vspace{0.1cm}

\item[$(v)$] \ \\
${\mathscr N}_{\epsilon}(V+\varphi,W+\psi)_K=0$.
\end{enumerate}
\end{prop}
\begin{proof}
Formula $(i)$ follows from Lemmas~\ref{brackets}, \ref{Lie deriv},
\ref{Half-diff}. Also, the vertical part of
$N_{\varepsilon}(A^h,B^h)_K$ is equal to
$$
\begin{array}{c}
{\mathcal V}{\mathscr N}_{\varepsilon}(A^h,B^h)_K
= -R(\pi_1(A),\pi_1(B))K +R(\pi_1(KA),\pi_1(KB))K \\[8pt]

\hspace{3.5cm}- {\mathscr K}_{\varepsilon}R(\pi_1(A),\pi_1(KB))K
-{\mathscr K}_{\varepsilon}R(\pi_1(KA),\pi_1(B))K.
\end{array}
$$
The part of ${\mathscr N}_{\varepsilon}(A^h,B^h)_K$ lying in ${\mathcal V}^{\ast}_Ê$ is the $1$-form whose value at every
vertical vector $W$ is
$$
\begin{array}{lll}
({\cal V}^{\ast}{\mathscr N}_{\varepsilon}(A^h,B^h)_K)(W)= \\[4pt]
-\displaystyle{\frac{1}{2}}[\pi_2(KA)(\pi_1(WB))+\pi_2(WB)(\pi_1(KA)) \\[6pt]
\hspace{2.5cm}-\pi_2(B)(\pi_1(({\mathscr P}_{\varepsilon}W)A))-\pi_2(({\mathscr P}_{\varepsilon}W)A)(\pi_1(B))] \\
                                                                           \\
+\displaystyle{\frac{1}{2}}[\pi_2(KB)(\pi_1(WA))+\pi_2(WA)(\pi_1(KB)) \\[6pt]
\hspace{2.5cm}-\pi_2(A)(\pi_1(({\mathscr P}_{\varepsilon}W)B))-\pi_2(({\mathscr P}_{\varepsilon}W)B)(\pi_1(A))]\\\\[6pt]

=-[<KA,WB>-<({\mathscr P}_{\varepsilon}W)A,B>] +[<KB,WA>-<({\mathscr P}_{\varepsilon}W)B,A>].
\end{array}
$$
Note also that
$$
<KA,WB>=<KW(A),B>=<{\mathscr P}_1W(A),B>.
$$
It follows that
$$
{\cal V}^{\ast}{\mathscr N}_{\varepsilon}(A^h,B^h)_J=-\omega^{\varepsilon}_{A,B}.
$$
This proves $(ii)$.

  To prove $(iii)$ take a section $a$ of $A(M)$ near the point $p$ such that $a(p)=V$
and $\nabla a|_p=0$. Let $\widetilde a$ be the vertical vector field defined by (\ref{eq tilde a}). Then it follows from
Lemmas~\ref{H-V brackets} and \ref{H-V Lie der&diff} that
$$
{\mathscr
N}_{\varepsilon}(A^h,V)_J=\frac{1}{2}N_{\varepsilon}(A^h,\widetilde
a)_K= (-({\mathscr P}_{\varepsilon}V)(A)+(K V)A)^{h}_K.
$$

  To prove $(iv)$,  take the vertical co-frame $\{\beta_{2t-1}^{\varepsilon},\beta_{2t}^{\varepsilon}\}$,
$t=1,...,n^2-n$, defined before the statement of Lemma~\ref{H-Vstar
Lie der}. Set
$\varphi=\sum\limits_{s=1}^{2(n^2-n)}\varphi_s^{\varepsilon}\beta_s^{\varepsilon}$,
$\varphi_s\in {\mathbb R}$. Let $E_1,..., E_{2n}$  be a basis of
$T_pM$ and $\xi_1,...,\xi_{2n}$ its dual basis.

We have $[A^h,\beta_s^{\varepsilon}]={\mathcal L}_{\pi_1(A^h)}\beta_s^{\varepsilon}$ since
$\beta_s^{\varepsilon}|{\mathcal H}=0$. Hence, by Lemma~\ref{H-Vstar Lie der}, $[A^h,\beta_s^{\varepsilon}]$ is the
horizontal form given by $Z^h\to -\beta_s^{\varepsilon}(R(\pi_1(A),Z)K)$. Thus
$$
[A^h,\beta_s^{\varepsilon}]_K=-\sum\limits_{k=1}^{2n}\beta_s^{\varepsilon}(R(\pi_1(A),E_k)K)(\xi_k)^h_K
$$
Similarly,
$$
\begin{array}{c}
[{\mathscr K}_{\varepsilon}A^h,{\mathscr K}_{\varepsilon}\beta_s^{\varepsilon}]_K=-\sum\limits_{k=1}^{2n}({\mathscr
K}_{\varepsilon}\beta_s^{\varepsilon})(R(\pi_1(KA),E_k)K)(\xi_k)^h_K\\[6pt]
= \sum\limits_{k=1}^{2n}\beta_s^{\varepsilon}({\mathscr K}_{\varepsilon}R(\pi_1(KA),E_k)K)(\xi_k)^h_K.
\end{array}
$$
$$
[{\mathscr
K}_{\varepsilon}A^h,\beta_s^{\varepsilon}]_K=-\sum\limits_{k=1}^{2n}\beta_s^{\varepsilon}(R(\pi_1(KA),E_k)K)(\xi_k)^h_K
$$
$$
[A^h,{\mathscr K}_{\varepsilon}\beta_s^{\varepsilon}]_K= \sum\limits_{k=1}^{2n}\beta_s^{\varepsilon}({\mathscr
K}_{\varepsilon}R(\pi_1(A),E_k)K)(\xi_k)^h_K
$$
Then
\begin{equation}\label{eq H-Vstar Nij}
\begin{array}{lll}
{\mathscr N}_{\varepsilon}(A^h,\varphi)_J=\sum\limits_{s=1}^{2(n^2-n)}\varphi_s^{\varepsilon} N_{\varepsilon}(A^h,\beta_s^{\varepsilon})_K=\\
                                                                                   \\
\sum\limits_{s=1}^{2(n^2-n)}\sum\limits_{k=1}^{2n}\varphi_s^{\varepsilon}\{[-\beta_s^{\varepsilon}(R(\pi_1(A),E_k)K)+
\beta_s^{\varepsilon}(K_{\varepsilon}R(\pi_1(KA),E_k)K)](\xi_k)^h_K \\
                                                                                    \\
+[\beta_s^{\varepsilon}(R(\pi_1(KA),E_k)K)-\beta_s^{\varepsilon}(K_{\varepsilon}R(\pi_1(A),E_k)K)](K\xi_k)^h_J\}.
\end{array}
\end{equation}
Moreover, note that
$$
<\xi_k,B>=\frac{1}{2}\xi_k(\pi_1(B)) \mbox { and } <K\xi_k,B>=-<\xi_k,KB>=-\frac{1}{2}\xi_k(\pi_1(KB)).
$$
Therefore
$$
\sum_{k=1}^{2n}<\xi_k,B>E_k=\frac{1}{2}\pi_1(B) \mbox { and }
\sum_{k=1}^{2n}<K\xi_k,B>E_k=-\frac{1}{2}\pi_1(KB).
$$
Now, $(iv)$ follows from (\ref{eq H-Vstar Nij}) and formula $(ii)$.

  Finally, identity $(v)$ follows from the fact that the generalized almost paracomplex structure
${\mathscr K}_{\varepsilon}$ on every fibre of ${\cal G}$ is induced by a paracomplex structure.
\end{proof}

\begin{cor}\label{non-integr}
The generalized almost paracomplex structures ${\mathscr K}_2, {\mathscr K}_3, {\mathscr K}_4$ are never integrable.
\end{cor}

\begin{proof}
Let $p\in M$ and take orthonormal bases $\{Q_1',...,Q_{2n}'\}$,
$\{Q_1'',...,Q_{2n}''\}$ of $E'_p$ and $E''_p$,  respectively. Let
$K$ be the product structure on $T_pM\oplus T_p^{\ast}M$ for which
$KQ'_{\alpha}=Q'_{n+\alpha}$ and $K'Q''_{\alpha}=Q''_{n+\alpha}$,
$\alpha=1,\dots,n$. Clearly, $K$ is a generalized paracomplex
structure on the vector space $T_pM$ compatible with the generalized
metric $E_p$. Define endomorphisms $S_{ij}'$ and $S_{ij}''$ by
(\ref{eq Sij}). Then $U'=S_{12}'+S_{n+1.n+2}'$ and
$U''=S_{12}''+S_{n+1.n+2}''$ are vertical tangent vectors of
${\mathcal G}$ at $K$. By Proposition~\ref{Nijenhuis} $(iii)$,
$N_2((Q_1'')^h,U'')=N_4((Q_1'')^h,U'')=2(Q_{n+2}'')^h$,
$N_3((Q_1')^h,U')=2(Q_{n+2}')^h$

\end{proof}

\medskip

Let ${\mathcal Z}={\mathcal Z}(TM,g)$ be the bundle over $M$ whose fibre at a point $p\in M$ consists of paracomplex
structures on $T_pM$ compatible with the metric $g$ (the usual reflector space of $(M,g)$). Denote by ${\mathcal N}^i$
the Nijenhuis tensor of the almost paracomplex structures ${\mathcal K}^i$ on ${\mathcal Z}$ defined in
Sec.~\ref{reflector}, i=1,2. To compute ${\mathcal N}^i$ we note first that Lemmas~\ref{brackets}  and \ref {H-V
brackets} have obvious analogs for the horizontal lifts to $T{\mathcal Z}$ of vector fields on $M$ and vertical vector
fields on ${\mathcal Z}$ defined by means of skew-symmetric endomorphisms of $TM$. Then a computation similar to that in
the proof of Proposition~\ref{Nijenhuis} gives the following.
\begin{prop}\label {Nij-Z}
Let $Q\in {\mathcal Z}$, $X,Y\in T_{\pi(Q)}M$, $V,W\in {\mathcal V}_Q$. Then:
$$
\begin{array}{c}
{\mathcal N}^i(X^h,Y^h)_Q=R(X,Y)Q+R(QX,QY)Q-{\mathcal K}^iR(QX,Y)Q-{\mathcal K}^iR(X,QY)Q;\\[6pt]
{\mathcal N}^i(X^h,V)_Q=[(-1)^i +1](QVX)^h_Q,\quad {\mathcal N}^i(V,W)=0.
\end{array}
$$

\end{prop}

\section{Integrability conditions in the case when the connection is determined by a generalized metric}\label{main case}

Let $D'=\nabla^{E'}$ be the connection on $E'=E$ determined by the
generalized metric $E$ (see Section 2.6). The image of this
connection under the isomorphism $pr_{TM}|E:E\to TM$ will be denoted
by $\nabla$. The connection $\nabla$ is metric with respect to the
metric $g$ and  has  skew-symmetric torsion:~ $g(T(X,Y),Z)=d\Theta
(X,Y,Z)$, $X,Y,Z\in TM$.

We define a connection $D''$ on $E''$  transferring  the connection
$\nabla$ from $TM$ to $E''$ by means of the isomorphism
$pr_{TM}|E'':E''\to TM$. Since this isomorphism is an isometry with
respect to the metrics $g''=-<.\,,.>|E''$ and $g$, we get a metric
connection on $E''$. As in the preceding section,  define a
connection $D$ on $TM\oplus T^{\ast}M$ setting $D=D'$ on $E'=E$ and
$D=D''$ on $E''$.

\smallskip

The connection $\nabla$ on $TM$ induces connections on the bundles obtained from $TM$ by algebraic operations like
$T^{\ast}M$, $TM\oplus T^{\ast}M$, etc. These connections will also be denoted by $\nabla$.

\smallskip

The connections $D$ and  $\nabla$ on $TM\oplus T^{\ast}M$ are
different in general. The connection $D$ preserves $E'$, while
$\nabla$ preserves $E'$ if and only if $\nabla\Theta=0$, a
condition, which is not always satisfied (take a closed form
$\Theta$, which is not parallel which respect the Levi-Civita
connection of $g$).

\smallskip

Recall that the reflector space  ${\mathcal Z}$ of $(M,g)$ is
considered as a subbundle of the bundle $A(TM)$ of
$g$-skew-symmetric endomorphisms of $TM$. The projections
$pr_{TM}|E':E'\to TM$ and $pr_{TM}|E'':E''\to TM$ yield an isometric
bundle-isomorphism $\chi:A(E')\oplus A(E'')\to A(TM)\oplus A(TM)$
sending the connection $D'\oplus D''$ to the connection
$\nabla\oplus \nabla$. The restriction of this map to the product
bundle ${\mathcal Z}(E')\times {\mathcal Z}(E'')\cong {\mathcal G}$
is the isomorphism onto the product bundle ${\mathcal Z}\times
{\mathcal Z}$ given by ${\mathcal G}\ni K\to (K_1,K_2)$, where
$K_1=(pr_{T}|E')\circ K \circ (pr_{T}|E')^{-1}$,
$K_2=(pr_{T}|E'')\circ K \circ (pr_{T}|E'')^{-1}$.

\smallskip

Suppose that $M$ is oriented and denote by ${\mathcal Z}_{\pm}$ the
bundle over $M$ whose fibre at a point $p\in M$ consists of
paracomplex structures on $T_pM$ compatible with the metric $g$ and
$\pm$ the orientation. Then the total spaces of the four product
bundles ${\mathcal Z}_{\pm}\times {\mathcal Z}_{\pm}$ are the
connected components of the manifold ${\mathcal Z}$. Thus ${\mathcal
G}$ has four connected component which the isomorphism
$\chi|{\mathcal G}:{\mathcal G}\to {\mathcal Z}\times {\mathcal Z}$
maps onto ${\mathcal Z}_{\pm}\times {\mathcal Z}_{\pm}$. The
connected component of ${\mathcal G}$ corresponding to ${\mathcal
Z}_{+}\times {\mathcal Z}_{+}$ will be denoted by ${\mathcal
G}_{++}$, and similar for the notations ${\mathcal G}_{+-}$,
${\mathcal G}_{-+}$, ${\mathcal G}_{--}$.

\begin{prop}\label {hor com zero}
$({\cal H}\oplus {\cal H}^{\ast}){\mathscr N}_{\varepsilon}(A^h,B^h)_K=0$ for every $K\in {\mathcal G}$ and every $A,B\in
TM\oplus T^{\ast}M$ if and only if $d\Theta=0$.
\end{prop}
\begin{proof}
  Let $K\in {\mathcal G}$ and let $S=(S',S'')$ be a section of ${\cal G}$ in a neighbourhood of the point $p=\pi(K)$
with the properties that $S(p)=K$ and $D S|_p=0$.

By Proposition ~\ref {Nijenhuis} $(i)$,  $({\cal H}\oplus {\cal
H}^{\ast}){\mathscr N}_{\varepsilon}(A^h,B^h)_J=0$ if and only if
the Nijenhuis tensor $N_S$ of the generalized almost paracomplex
structure $S$ on $M$ vanishes at the  point $p$. Let $S_1$ and $S_2$
be the almost paracomplex structures on $M$ determined by $S$,
$$
S_1=(\pi_1|E')\circ S'\circ (\pi_1|E')^{-1},\quad S_2=(\pi_1|E'')\circ S''\circ (\pi_1|E'')^{-1}.
$$
These structures are compatible with the metric $g$ and we denote their fundamental $2$-forms  by $\Omega_1$ and
$\Omega_2$, respectively:
$$
\Omega_1(X,Y)=g(X,S_1Y),\quad \Omega_2(X,Y)=g(X,S_2Y),\quad X,Y\in TM.
$$

 Denote by $P$ the generalized paracomplex structure on $M$ with the
block-matrix
$$P=\frac{1}{2}
  \left(
  \begin{array}{cc}
    S_{1}+S_{2} & -\Omega_{1}^{-1}+\Omega_{2}^{-1} \\
    -\Omega_{1}+\Omega_{2}& -S_{1}^{\ast}-S_{2}^{\ast} \\
  \end{array}
 \right)
$$
By Proposition~\ref{K-K1,2}, the generalized paracomplex structure $S$ is the $B$-transform of $P$ by means of the form
$\Theta$:
$$
S=e^{\Theta}Pe^{-\Theta}.
$$
Let $N_P$ be the Nijensuis tensor of the generalized almost
paracomplex structure $P$. Set
\begin{equation}\label{notat}
A=X+\alpha,\quad B=Y+\beta,\quad PA=\widehat{X}+\widehat{\alpha},\quad PB=\widehat{Y}+\widehat{\beta}
\end{equation}
where $X,Y,\widehat{X},\widehat{Y}\in TM$ and $\alpha,\beta,\widehat{\alpha},\widehat{\beta}\in T^{\ast}M$. Then, by
Proposition~\ref{Courant-B-transf} and the fact that $e^{-\Theta}|T^{\ast}M=Id$,
$$
\begin{array}{c}
N_S(e^{\Theta}A,e^{\Theta}B)=e^{\Theta}N_P(A,B)+\imath_{Y}\imath_{X}d\Theta
+\imath_{\widehat{Y}}\imath_{\widehat{X}}d\Theta\\[8pt]
\hspace{7cm}-e^{\Theta}P(\imath_{Y}\imath_{\widehat{X}}d\Theta +\imath_{\widehat{Y}}\imath_{X}d\Theta).
\end{array}
$$
It follows that $({\cal H}\oplus {\cal H}^{\ast})N_{\varepsilon}(A^h,B^h)_K=0$ for every $A,B\in TM\oplus T^{\ast}M$ if
and only if at the point $p=\pi(K)$
\begin{equation}\label{cond}
N_P(A,B)=-\imath_{Y}\imath_{X}d\Theta -\imath_{\widehat{Y}}\imath_{\widehat{X}}d\Theta
+P(\imath_{Y}\imath_{\widehat{X}}d\Theta+\imath_{\widehat{Y}}\imath_{X}d\Theta).
\end{equation}
 We have
$$
\nabla S_1=(\pi_1|E')\circ (D S')\circ (\pi_1|E')^{-1}
$$
since the connection $\nabla$ on $TM$ is obtained from the
connection $D|E'=\nabla^{E'}$ by means of the isomorphism
$\pi_1|E':E'\to TM$. In particular $\nabla S_1|_p=0$. Similarly,
$\nabla S_2|_p=0$. Then $\nabla S_k^{\ast}|_p=0$, $k=1,2$, and
$\nabla\Omega_k|_p=-\nabla(g\circ S_k)|_p=0$,
$\nabla\Omega_k^{-1}=\nabla(S_k\circ g^{-1})|_p=0$ since $\nabla
g=\nabla g^{-1}=0$. It follows that $\nabla P|_p=0$. Extend $X$ and
$\alpha$ to a vector field $X$ and a $1$-form $\alpha$ on $M$ such
that $\nabla X|_p=0$ and $\nabla\alpha|_p=0$; similarly for $Y$ and
$\beta$. In this way, we obtain sections $A=X+\alpha$ and
$B=Y+\beta$ of $TM\oplus T^{\ast}M$ such that  $\nabla A|_p=\nabla
B|_p=0$ and $\nabla PA|_p=\nabla PB|_p=0$

In order to compute $N_P(A,B)$, we need the following simple
observation. Let $Z$ be a vector field and $\omega$ a $1$-form on
$M$ such that $\nabla Z|_p=0$ and $\nabla\omega|_p=0$. Then, for
every $Z'\in T_pM$,
$$
\begin{array}{c}
({\cal L}_{Z}\omega)(Z')_p
=(\nabla_{Z}\omega)(Z')_p+\omega(T(Z,Z')=\omega(T(Z,Z'))\\[8pt]
(d\,\imath_{Z}\omega)(Z')_p=Z'_p(\omega(Z))=(\nabla_{Z'}\omega)(Z)_p=0
\end{array}
$$
where $T(Z,Z')$ is the torsion tensor of the connection $\nabla$.
For $Z\in TM$, let $\imath_{Z}T:TM\to TM$ be the map $Z'\to
T(Z,Z')$. Then, under the notation of (\ref{notat}),
\begin{equation}\label{NP}
\begin{array}{c}
N_P(A,B)=-T(X,Y)-T(\widehat{X},\widehat{Y})-\alpha(\imath_{Y}T)+\beta(\imath_{X}T)
-\widehat{\alpha}(\imath_{\widehat{Y}}T)+\widehat{\beta}(\imath_{\widehat{X}}T)\\[8pt]
+P[T(\widehat{X},Y)+T(X,\widehat{Y})+\widehat{\alpha}(\imath_{Y}T)-\beta(\imath_{\widehat{X}}T)
+\alpha(\imath_{\widehat{Y}}T)-\widehat{\beta}(\imath_{X}T)].
\end{array}
\end{equation}

If $\alpha=g(X')$ for some (unique)  $X'\in TM$,
$$
\alpha(\imath_{Y}T)=g(T(X',Y))=\imath_{Y}\imath_{X'}d\Theta,\quad Y\in TM,
$$
and
$$
P(\alpha\circ\imath_{Y}T)=-\displaystyle{\frac{1}{2}}[(\Omega_1^{-1}-\Omega_2^{-1})\imath_{Y}\imath_{X'}d\Theta
+(S_1^{\ast}+S_2^{\ast})\imath_{Y}\imath_{X'}d\Theta]
$$

Moreover, $g(T(X,Y))=\imath_{Y}\imath_{X}d\Theta$  for every $X,Y\in TM$, hence
$$
P(T(X,Y))=-\displaystyle{\frac{1}{2}}[(\Omega_1^{-1}+\Omega_2^{-1})\imath_{Y}\imath_{X}d\Theta
+(S_1^{\ast}-S_2^{\ast})\imath_{Y}\imath_{X}d\Theta].
$$
Note also that
$$
\begin{array}{c}
\widehat{X}=\displaystyle{\frac{1}{2}}[S_1(X+X')+S_2(X-X')],\\[8pt]
\widehat{\alpha}=\displaystyle{\frac{1}{2}}[g(S_1(X+X')-S_2(X-X'))].
\end{array}
$$

Now, suppose that
$$
({\cal H}\oplus {\cal H}^{\ast})N_{\varepsilon}(A^h,B^h)_K=0, \quad A,B\in TM\oplus T^{\ast}M.
$$
Then, by (\ref{cond}),
\begin{equation}\label{cond-2-forms}
N_P(g(X),g(Y))=\displaystyle{\frac{1}{4}\imath_{(S_1X-S_2X)}\imath_{(S_1Y-S_2Y)}d\Theta}.
\end{equation}
Therefore the tangential component of $N_P(g(X),g(Y))$ vanishes. Hence by (\ref{NP})
$$
\begin{array}{c}
-T(S_1X-S_2X,S_Y-S_2Y)\\[6pt]
+(\Omega^{-1}_1-\Omega^{-1}_2)\imath_{(S_1X-S_2X)}\imath_{Y}d\Theta
-(\Omega^{-1}_1-\Omega^{-1}_2)\imath_{(S_1Y-S_2Y)}\imath_{X}d\Theta=0.
\end{array}
$$
Applying the map $g$ to both sides of the latter identity, we obtain
by means of the identities $g\circ\Omega_k^{-1}=S_k^{\ast}$ that for
every $X,Y,Z\in T_pM$
\begin{equation}\label{ff1}
\begin{array}{c}
d\Theta(S_1X-S_2X,Y,S_1Z-S_2Z)+d\Theta(X,S_1Y-S_2Y,S_1Z-S_2Z)\\[8pt]
=-d\Theta(S_1X-S_2X,S_1Y-S_2Y,Z).
\end{array}
\end{equation}
Applying (\ref{ff1}) for the generalized almost paracomplex
structure determined by the paracomplex structures $(-S_1,S_2)$ on
$T_pM$ and comparing the obtained identity with (\ref{ff1}), we see
that
\begin{equation}\label{dTheta-1}
\begin{array}{c}
d\Theta(S_1X,Y,S_1Z)+d\Theta(S_2X,Y,S_2Z)+d\Theta(X,S_1Y,S_1Z)+d\Theta(X,S_2Y,S_2Z)\\[8pt]
=-d\Theta(S_1X,S_1Y,Z)-d\Theta(S_2X,S_2Y,Z).
\end{array}
\end{equation}
Computing the co-tangential component of $N_P(g(X),g(Y))$ by means of (\ref{NP}), we get from (\ref{cond-2-forms})
$$
\begin{array}{c}
-\displaystyle{\imath_{(S_1Y-S_2Y)}\imath_{(S_1X+S_2X)}d\Theta +\imath_{(S_1X-S_2X)}\imath_{(S_1Y+S_2Y)}d\Theta}\\[6pt]
+\displaystyle{(S^{\ast}_1+S^{\ast}_2)\big(\imath_{(S_1X-S_2X)}\imath_{Y}d\Theta-\imath_{(S_1Y-S_2Y)}\imath_{X}d\Theta\big)}\\[6pt]
=\displaystyle{\imath_{(S_1X-S_2X)}\imath_{(S_1Y-S_2Y)}d\Theta}.
\end{array}
$$
Taking the value of both sides at a vector $Z\in T_pM$, then
applying the obtained identity for the generalized almost
paracomplex structure determined by $(-S_1,S_2)$, we obtain
\begin{equation}\label{dTheta-2}
\begin{array}{c}
-d\Theta(S_1X,Y,S_1Z)+d\Theta(S_2X,Y,S_2Z)-d\Theta(X,S_1Y,S_1Z)+d\Theta(X,S_2Y,S_2Z)\\[8pt]
=d\Theta(S_1X,S_1Y,Z)-3d\Theta(S_2X,S_2Y,Z).
\end{array}
\end{equation}
It follows from (\ref{dTheta-1}) and (\ref{dTheta-2}) that
$$
d\Theta(S_2X,Y,S_2Z)+d\Theta(X,S_2Y,S_2Z)=-2d\Theta(S_2X,S_2Y,Z)
$$
Hence
\begin{equation}\label{eq1}
-2d\Theta(X,Y,Z)=d\Theta(X,S_2Y,S_2Z)+d\Theta(S_2X,Y,S_2Z).
\end{equation}
This implies
$$
-[d\Theta(X,S_2Y,S_2Z)+d\Theta(S_2X,Y,S_2Z)]=d\Theta(X,Y,Z)+d\Theta(S_2X,S_2Y,Z)
$$
The latter identity and (\ref{eq1}) give
\begin{equation}\label{eq2}
d\Theta(X,Y,Z)=d\Theta(S_2X,S_2Y,Z).
\end{equation}
Therefore
$$
d\Theta(X,S_2Y,S_2Z)=d\Theta(S_2X,Y,S_2Z).
$$
Then, by (\ref{eq1}),
\begin{equation}\label{eq3}
d\Theta(X,Y,Z)=-d\Theta(X,S_2Y,S_2Z)
\end{equation}
Let $e_1,e_2,...,e_{2n-1},e_{2n}$ be an orthonormal basis of $T_pM$
with $g(e_{2i-1},e_{2j-1})=-g(e_{2i},e_{2j})=\delta_{ij}$,
$g(e_{2i-1},e_{2j-1})=0$, $i,j=1,...,n$. For every fixed $i$ and
$j$, there is a paracomplex structure $K$ on $T_pM$ compatible with
the metric such that $Ke_{2i-1}=e_{2j}$. Identity (\ref{eq2}) with
$S_2=K$, $X=e_{2i-1}$, $Y=e_{2j}$, and arbitrary $Z\in T_pM$ gives
\begin{equation}\label{odd-even}
d\Theta(e_{2i-1},e_{2j},Z)=0,\quad Z\in T_pM.
\end{equation}
Now, fix $j$ and take a compatible paracomplex structure $K$ on
$T_pM$ such that $K_{2j-1}=e_{2j}$, Then by (\ref{eq3}) and
(\ref{odd-even})
$$
d\Theta(e_{2i-1},e_{2j-1},Z)=-d\Theta(e_{2i-1},e_{2j},KZ)=0.
$$
The latter identity and (\ref{eq3}) imply
$$
d\Theta(e_{2i},e_{2j},Z)=d\Theta(e_{2i},e_{2j-1},KZ)=0,
$$
Hence $d\Theta(e_{\alpha},e_{\beta},Z)=0$ for every
$\alpha,\beta=1,...,2n$ and every $Z\in T_pM$, Thus, $d\Theta=0$.

Conversely, if $d\Theta=0$, then $T=0$ and we have $N_K=0$ by
(\ref{NP}).  Thus, the condition (\ref{cond}) is trivially
satisfied. Therefore
$$
({\cal H}\oplus {\cal H}^{\ast})N_{\varepsilon}(A^h,B^h)_K=0, \quad A,B\in TM\oplus T^{\ast}M.
$$
\end{proof}

Suppose that $M$ is oriented and $dim\,M=4k$. Then the above proof still holds true if we, instead of ${\mathcal G}$,
consider a connected component of it. Indeed, the almost paracomplex structures $S_1$ and $-S_1$ above induce the same
orientation.  Moreover, if $e_1,e_2,...,e_{2n-1},e_{2n}$ is an orthonormal basis of $T_pM$ with
$g(e_{2i-1},e_{2j-1})=-g(e_{2i},e_{2j})=\delta_{ij}$, $g(e_{2i-1},e_{2j-1})=0$, $i,j=1,...,n$, and if we fix two indexes
$i$, $j$,  the paracomplex structure $K$ on $T_pM$ compatible with the metric with the property that $Ke_{2i-1}=e_{2j}$
used at the end of the proof can be chosen to induce the given or the opposite orientation of $M$. Thus we have the
following.

\begin{prop}\label{dim=4k} If $M$ is oriented and $dim\,M=4k$, then
$$({\cal H}\oplus {\cal H}^{\ast}){\mathscr N}_{\varepsilon}(A^h,B^h)_K=0$$
for every $K$ in a connected component of ${\mathcal G}$ and every $A,B\in TM\oplus T^{\ast}M$ if and only if
$d\Theta=0$.
\end{prop}

\smallskip

Recall that if $R$ is the curvature tensor of the Levi-Civita connection of $(M,g)$,  the curvature operator ${\mathcal
R}$ is the self-adjoint endomorphism of $\Lambda ^2TM$ defined by
$$
g({\mathcal R}(X\land Y),Z\land T)=g(R(X,Y)Z,T),\quad X,Y,Z,T\in TM.
$$
The metric on $\Lambda^2TM$ used in the left-hand side of the latter identity is defined by
$$
g(X_1\wedge X_2,X_3\wedge X_4)=g(X_1,X_3)g(X_2,X_4)-g(X_1,X_4)g(X_3,X_4).
$$
As is well-known, the curvature operator decomposes as (see, for example, \cite [Section 1 G, H]{Besse})
\begin{equation}\label{curv decom}
{\mathcal R}=\frac{s}{n(n-1)}Id + {\mathcal B}+ {\mathcal W}
\end{equation}
where $s$ is the scalar curvature of the manifold $(M,g)$ and ${\mathcal B}$, ${\mathcal W}$ correspond to its traceless
Ricci tensor and Weyl conformal tensor, respectively.

If $\rho:TM\to TM$ is the Ricci operator, $g(\rho(X),Y)=Ricci(X,Y)$, the operator ${\mathcal B}$ is given by
\begin{equation}\label{B}
{\mathcal B}(X\wedge Y)=\frac{1}{n-2}[\rho(X)\wedge Y + X\wedge\rho(Y) -\frac{2s}{n} X\wedge Y], \quad X,Y\in TM.
\end{equation}
Thus, a Riemannian manifold  is Einstein exactly when ${\cal B}=0$; it is conformally flat when ${\mathcal W}=0$.

\smallskip

If the dimension of $M$ is four, the Hodge star operator defines an involution $\ast$ of $\Lambda^2TM$ and we have the
orthogonal decomposition
$$
\Lambda^2TM=\Lambda^2_{-}TM\oplus\Lambda^2_{+}TM
$$
where $\Lambda^2_{\pm}TM$ are the subbundles of $\Lambda^2TM$ corresponding to the $(\pm 1)$-eigenvalues of the operator
$\ast$. Accordingly, the operator ${\mathcal W}$ has an extra decomposition ${\mathcal W}={\mathcal W}_{+}+{\mathcal
W}_{-}$ where ${\mathcal W}_{\pm}={\mathcal W}$ on $\Lambda^2_{\pm}TM$ and ${\mathcal W}_{\pm}=0$ on $\Lambda^2_{\mp}TM$.
The operator ${\mathcal B}$ does not have a decomposition of this type since it maps $\Lambda^2_{\pm}TM$ into
$\Lambda^2_{\mp}TM$. Reversing the orientation of $M$ interchanges the roles of $\Lambda^2_{+}TM$ and $\Lambda^2_{-}TM$,
hence the roles of ${\mathcal W}_{+}$ and ${\mathcal W}_{-}$.

Recall also that  $(M,g)$ is  called self-dual (anti-self-dual), if ${\cal W}_{-}=0$ (resp. ${\cal W}_{+}=0$). As is
well-known, the self-duality (anti-self-duality) is an interesting curvature condition which plays an important  role in
twistor theory.

\subsection {The case when the base manifold is four-dimension}

\smallskip

\noindent \\
Henceforth we assume that $M$ is oriented and of dimension $4$.

\smallskip

According to Propositions~\ref{Nijenhuis} and \ref{dim=4k}, the restriction of the generalized almost paracomplex
structure ${\mathscr K}_1$ to a connected component $\widetilde{\mathcal G}$ of ${\mathcal G}$ is integrable if and only
if $d\Theta=0$ and for every $p\in M$, $A,B\in T_pM$, and for every generalized paracomplex structure
$K\in\widetilde{\mathcal G}$ on $T_pM$
$$
\begin{array}{c}
R(\pi_1(A),\pi_1(B))K +R(\pi_1(KA),\pi_1(KB))K \\[8pt]
-{\mathscr P}_1R(\pi_1(KA),\pi_1(B))K - {\mathscr P}_1R(\pi_1(A),\pi_1(KB))K =0
\end{array}
$$
where $R$ is the curvature tensor of the connection $D$ on the bundle $A(E')\oplus A(E'')$. If $(K_1,K_2)$ are the
paracomplex structures on $T_pM$ determined by $K$, the latter identity is equivalent to the identities
\begin{equation}\label{I and II comp}
\begin{array}{c}
R(\pi_1(A),\pi_1(B))K_r +R(\pi_1(KA),\pi_1(KB))K_r \\[8pt]
-K_r\circ R(\pi_1(KA),\pi_1(B))K_r- K_r\circ R(\pi_1(A),\pi_1(KB))K_r =0, \quad r=1,2,
\end{array}
\end{equation}
where $R$ is the curvature tensor on the bundle $A(TM)$ of skew-symmetric endomorphism of $TM$ induced by the connection
$\nabla$.

Assume that $d\Theta=0$. Then $\nabla$ is the Levi-Civita connection
of the Riemannian manifold $(M,g)$. Every $A\in E_p'$ is of the form
$A=X+g(X)+\Theta(X)$ for some (unique) $X\in T_pM$ and
$KA=K_1X+g(K_1X)+\Theta(K_1X)$. Similarly, if $B\in E''_p$, then
$B=Y-g(Y)+\Theta(Y)$, $Y\in T_pM$, and
$KB=K_2Y-g(K_2Y)+\Theta(K_2Y)$.  It follows that the identity
(\ref{I and II comp}) is equivalent to the condition that for every
$X,Y,Z,U\in T_pM$ and every paracomplex structures $(K_1,K_2)$ on
$T_pM$ corresponding to a generalized paracomplex structure $K$ in
$\widetilde{\mathcal G}$
\begin{equation}\label{Jklr}
\begin{array}{c}
g({\mathcal R}(X\wedge Y+K_jX\wedge K_lY), Z\wedge U+K_rZ\wedge K_rU)\\[6pt]
\hspace{3.3cm}=-g({\mathcal R}(K_jX\wedge Y+X\wedge K_lY),K_rZ\wedge U+Z\wedge K_rU),\\[6pt]
 j,l,r=1,2.
\end{array}
\end{equation}
The paracomplex structures $(K_1,K_2)$ in the latter identities are compatible with the metric $g$ and, moreover, they
induce the orientation of $T_pM$ if we consider the connected component $\widetilde{\mathcal G}={\mathcal G}_{++}$, while
$(K_1,K_2)$ induce the opposite orientation in the case $\widetilde{\mathcal G}={\mathcal G}_{--}$. If
$\widetilde{\mathcal G}={\mathcal G}_{+-}$, the paracomplex structure $K_1$ induces the given orientation of $T_pM$ and
$K_2$ yields the opposite one, and vice versa if $\widetilde{\mathcal G}={\mathcal G}_{-+}$.

\smallskip

According to Proposition~\ref{Nij-Z}, identity (\ref{Jklr}) for
$j=l=r$  coincides with the integrability condition for the almost
paracomplex structure ${\mathcal K}_1$ on the positive or negative
reflector space of $(M,g)$. By \cite[Theorem 5]{JR}, this
integrability condition is the anti-self-duality of $(M,g)$ in the
case of positive reflector spaces and is the self-duality of $(M,g)$
when considering the negative reflector space.

\begin{theorem}\label{++}
The restriction of the generalized paracomplex structure ${\mathscr K}_1$ to ${\mathcal G}_{++}$ (resp. ${\mathcal
G}_{--}$)  is integrable if and only if $(M,g)$ is anti-self-dual (resp. self-dual) and Ricci flat.
\end{theorem}

\begin{proof} Let $E_1,..., E_4$ be an oriented orthonormal basis of a tangent space $T_pM$ with $||E_1||^2=||E_2||^2=1$,
$||E_3||^2=||E_4||^2=-1$.  It is convenient to set $E_{ab}=E_a\wedge E_b$ and $\rho_{ab}=Ricci(E_a,E_b)$, $a,b=1,...,n$.

Suppose that the structure ${\mathscr K}_1| {\mathcal G}_{++}$ is integrable. Then ${\mathcal W}_{+}=0$ as we have
remarked above. Let $K_1$ and $K_2$ be paracomplex structures on $T_pM$ for which $K_1E_{1}=E_{3}$, $K_1E_{2}=E_{4}$ and
$K_2E_{1}=E_{4}$, $K_2E_{2}=-E_{3}$. Identity (\ref{Jklr}) with $j=l=1$, $r=2$, and $(X,Y,Z,U)=(E_1,E_2,E_3,E_4)$ gives
\begin{equation}\label{k=l=1,r=2}
g({\mathcal R}(E_{12}+E_{34}),E_{12}+E_{34})-g({\mathcal R}(E_{14}-E_{23}),E_{13}-E_{42})=0.
\end{equation}
 Since ${\mathcal W}_{+}=0$ and $E_{12}+E_{34},E_{14}-E_{23}\in\Lambda^2_{+}T_pM$, we have
$$
{\mathcal W}(E_{12}+E_{34})= {\mathcal W}(E_{14}-E_{23})=0.
$$
Hence
$$
\begin{array}{c}
g({\mathcal R}(E_{12}+E_{34}),E_{12}+E_{34})-g({\mathcal
R}(E_{14}-E_{23}),E_{13}-E_{42})\\[6pt]
=\displaystyle{\frac{s}{6}}+g({\mathcal B}(E_{12}+E_{34}),E_{12}+E_{34})+g({\mathcal B}(E_{14}-E_{23}),E_{13}-E_{42}).
\end{array}
$$
By (\ref{B})
$$
\begin{array}{c}
g({\mathcal
B}(E_{12}+E_{34}),E_{12}+E_{34})=\displaystyle{\frac{1}{2}}[\rho_{11}+\rho_{22}-\rho_{33}-\rho_{44}-s]=0,\\[8pt]
g({\mathcal B}(E_{14}-E_{23}),E_{13}-E_{42})=0.
\end{array}
$$
Then by (\ref{k=l=1,r=2}) $s=0$.

In order to show that ${\mathcal B}=0$, we apply identity
(\ref{Jklr}) with $j=1$, $l=2$ and  take $K_1$, $K_2$ to be the
paracomplex structures introduced above. Summing the identities
corresponding to $X=Y=E_2$ and $X=Y=E_3$, we get
$$
g({\mathcal R}(E_{12}-E_{34}),Z\wedge U+K_rZ\wedge K_rU)=g({\mathcal R}(E_{13}+E_{42}),K_rZ\wedge U+Z\wedge K_rU).
$$
Summing the identities corresponding to $X=Y=E_1$ and $X=Y=E_4$, we
obtain
$$
g({\mathcal R}(E_{12}-E_{34}),Z\wedge U+K_rZ\wedge K_rU)=-g({\mathcal R}(E_{13}+E_{42}),K_rZ\wedge U+Z\wedge K_rU).
$$
Thus,
\begin{equation}\label{aux-1}
g({\mathcal R}(E_{12}-E_{34}),Z\wedge U+K_rZ\wedge K_rU)=0=g({\mathcal R}(E_{13}+E_{42}),K_rZ\wedge U+Z\wedge K_rU).
\end{equation}
We also apply identity (\ref{Jklr}) for $K_j=K_2$, $K_l=-K_2$. Subtracting the identities corresponding to $X=Y=E_1$ and
$X=Y=E_2$, we get
\begin{equation}\label{aux-2}
g({\mathcal R}(E_{14}+E_{23}),K_rZ\wedge U+Z\wedge K_rU)=0.
\end{equation}
Every $2$-vector of the form $Z\wedge U+K_rZ\wedge K_rU$ lies in
$\Lambda^{2}_{+}T_pM$. Indeed, this $2$-vector is orthogonal to
$\Lambda^2_{-}T_pM$ since $K_r$ commutes with every endomorphism of
$T_pM$ corresponding to a $2$-vector in $\Lambda^2_{-}T_pM$ via
(\ref{Sa}). Therefore ${\mathcal W}(Z\wedge U+K_rZ\wedge
K_rU)={\mathcal W}( K_rZ\wedge U+Z\wedge K_rU)=0$. Now,  set $r=1$,
$(Z,U)=(E_1,E_2)$ in the first identity in (\ref{aux-1}), $r=2$,
$(Z,U)=(E_4,E_3)$ in the second one, and $r=1$, $(Z,U)=(E_3,E_4)$ in
(\ref{aux-2}). This gives
$$
\begin{array}{c}
g({\mathcal B}(E_{12}-E_{34}),E_{12}+E_{34})=g({\mathcal B}(E_{13}+E_{42}),E_{13}-E_{42})\\[6pt]
=g({\mathcal B}(E_{14}+E_{23}),E_{14}-E_{23})=0.
\end{array}
$$
By (\ref {B}), it follows
$$
\rho_{11}+\rho_{22}+\rho_{33}+\rho_{44}=0,\quad \rho_{11}-\rho_{22}-\rho_{33}+\rho_{44}=0,\quad
\rho_{11}-\rho_{22}+\rho_{33}-\rho_{44}=0.
$$
These identities, together with
$\rho_{11}+\rho_{22}-\rho_{33}-\rho_{44}=s=0$, imply
$\rho_{11}=\rho_{22}=\rho_{33}=\rho_{44}=0$.  Moreover, the first
identity in (\ref{aux-1}) for $r=1$ and $(Z,U)=(E_1,E_4)$ reads as
$$
g({\mathcal B}(E_{13}-E_{42}),E_{12}+E_{34})=0.
$$
This gives $\rho_{24}+\rho_{13}=0$. Similarly, it follows from
(\ref{aux-2}) with $r=1$ and $(Z,U)=(E_2,E_3)$ that
$\rho_{24}-\rho_{13}=0$. Hence $\rho_{13}=\rho_{24}=0$. Replacing
the basis $(E_1,E_2,E_3,E_4)$ by $(E_1,E_2,E_4,-E_3)$, we see that
$\rho_{14}=\rho_{23}=0$. The second identity of (\ref{aux-1}) with
$r=1$, $(Z,U)=(E_3,E_4)$ gives $\rho_{12}+\rho_{34}=0$; identity
(\ref{aux-2}) with $r=2$ and $(Z,U)=(E_3,E_4)$ implies
$\rho_{12}-\rho_{34}=0$. Thus, $\rho_{12}=\rho_{34}=0$. This shows
that $\rho_{ij}=0$ for every $i,j=1,...,4$, i.e.,  $Ricci=0$.

Conversely, if $s=0$ and ${\mathcal B}={\mathcal W}_{+}=0$, identity (\ref{Jklr}) is trivially satisfied since for every
$X,Y\in T_pM$ and every paracomplex structure $K$ on $T_pM$ compatible with the metric and orientation, the $2$-vector
$X\wedge Y+KX\wedge KY$ lies in $\Lambda^2_{+}T_pM$, so ${\mathcal R}(X\wedge Y+KX\wedge KY)=0$.

\end{proof}

\begin{theorem}\label{+-}
The restriction of the generalized paracomplex structure ${\mathscr K}_1$ to ${\mathcal G}_{++}$ (resp. ${\mathcal
G}_{--}$)  is integrable if and only if $(M,g)$ is of constant sectional curvature.
\end{theorem}

\begin{proof}
Suppose that ${\mathscr K}_1|{\mathcal G}_{+-}$ is integrable. Then, by the remarks preceding the statement of
Theorem~\ref{++},  $(M,g)$ is both anti-self-dual and self-dual, hence ${\mathcal W}=0$. Next, note that  that, by
(\ref{B}), the condition ${\mathcal B}=0$ is equivalent to
$$
\rho_{aa}||E_a||^2+\rho_{bb}||E_b||^2-\frac{s}{2}=0, \quad \rho_{ab}=0,\quad a\neq b,\quad a,b=1,...,4.
$$
In order to prove that ${\mathcal B}=0$, let  $E_1,..., E_4$ be an
oriented orthonormal basis of a tangent space $T_pM$ with
$||E_1||^2=||E_2||^2=1$, $||E_3||^2=||E_4||^2=-1$. Apply identity
(\ref{Jklr}) with $j=1$ and $l=2$, and  take for $K_1$ and $K_2$ the
paracomplex structures on $T_pM$ for which $K_1E_{1}=E_{3}$,
$K_1E_{4}=E_{2}$ and $K_2E_{1}=E_{4}$, $K_2E_{2}=E_{3}$. Add the
identities corresponding to: $(1)$ $(X,Y)=(E_1,E_4)$ and
$(X,Y)=(E_2,E_3)$, $(2)$ $(X,Y)=(E_1,E_4)$ and $(X,Y)=(E_3,E_2)$. In
this way, we get
\begin{equation}\label{eq1}
g({\mathcal R}(E_{14}+E_{23}),Z\wedge U+K_rZ\wedge K_rU)-g({\mathcal R}(E_{13}-E_{42}),Z\wedge U+K_rZ\wedge K_rU)=0.
\end{equation}
\begin{equation}\label{eq2}
g({\mathcal R}(E_{14}-E_{23}),Z\wedge U+K_rZ\wedge K_rU)+g({\mathcal R}(E_{12}+E_{34}),K_rZ\wedge U + Z\wedge K_rU)=0
\end{equation}
Put $r=1$ in (\ref{eq1})  and set $(Z,U)=(E_1,E_4)$. This gives
\begin{equation}\label{+-1}
g({\mathcal R}(E_{14}+E_{23}),E_{14}-E_{23})-g({\mathcal R}(E_{13}-E_{42}),E_{14}-E_{23})=0.
\end{equation}
Setting $r=2$ and  $(Z,U)=(E_1,E_3)$ in (\ref{eq1}), we obtain
\begin{equation}\label{+-2}
g({\mathcal R}(E_{14}+E_{23}),E_{13}+E_{42})-g({\mathcal R}(E_{13}-E_{42}),E_{13}+E_{42})=0.
\end{equation}
Now,  apply (\ref{eq2}) for the paracomplex structure $K_2$
determined by $K_2E_1=E_3$, $KE_2=-E_4$. Setting $(Z,U)=(E_1,E_4)$
gives
\begin{equation}\label{mix}
g({\mathcal R}(E_{14}-E_{23}),E_{14}+E_{23})-g({\mathcal R}(E_{12}+E_{34}),E_{12}-E_{34})=0
\end{equation}

The second term in (\ref{+-1}) and the first one in (\ref{+-2}) vanish since ${\mathcal W}=0$. It follows that
$$
g({\mathcal B}(E_{14}+E_{23}),E_{14}-E_{23})=0,\quad g({\mathcal B}(E_{13}-E_{42}),E_{13}+E_{42})=0.
$$
The first of this identities and (\ref{mix}) imply
$$
g({\mathcal B}(E_{12}+E_{34}),E_{12}-E_{34})=0
$$
Then, by (\ref{B}),
$$
\rho_{11}-\rho_{22}+\rho_{33}-\rho_{44}=0,\quad \rho_{11}-\rho_{22}-\rho_{33}+\rho_{44}=0,\quad
\rho_{11}+\rho_{22}+\rho_{33}+\rho_{44}=0.
$$
This, together with the identity
$\rho_{11}+\rho_{22}-\rho_{33}-\rho_{44}=s$, implies
$$
\rho_{11}=\rho_{22}=\frac{s}{4},\quad \rho_{33}=\rho_{44}=-\frac{s}{4}.
$$
Next, we use (\ref{Jklr}) with $j=1$, $l=2$ and $K_1$, $K_2$ defined
by $K_1E_1=E_3$, $K_1E_2=E_4$, $K_2E_{1}=E_{4}$, $K_2E_{2}=E_{3}$.
Setting $(X,Y)=(E_2,E_3)$, we get
$$
g({\mathcal R}(E_{23}),Z\wedge U+K_rZ\wedge K_rU)=0.
$$
Set $(Z,U)=(E_1,E_2)$ in the latter identity and apply it twice -
with $r=1$ and $K_1$ defined by $K_1E_1=E_3$, $K_1E_2=E_4$, as well
as with $r=2$ and $K_2$ defined by $K_2E_1=E_3$, $K_2E_2=-E_4$. Then
 set $(Z,U)=(E_1,E_3)$ and take $K_1$ and $K_2$ determined by
$K_1E_1=E_4$, $K_1E_2=-E_3$ and $K_2E_1=E_4$, $K_2E_2=E_3$. Thus, we
obtain
$$
\begin{array}{c}
g({\mathcal R}(E_{23},E_{12}+E_{34})=g({\mathcal R}(E_{23},E_{12}-E_{34})=0,\\[6pt]
g({\mathcal R}(E_{23},E_{13}-E_{42})=g({\mathcal R}(E_{23},E_{13}+E_{42})=0.
\end{array}
$$
These identities imply
$$
\rho_{24}-\rho_{13}=-\rho_{24}-\rho_{13}=0,\quad -\rho_{12}+\rho_{34}=-\rho_{12}+\rho_{34}=0.
$$
Hence $\rho_{13}=\rho_{24}=\rho_{12}=\rho_{14}=0$. Replacing the
basis $(E_1,E_2,E_3,E_4)$ by $(E_2,E_1,E_3,-E_4)$, we see that
$\rho_{23}=\rho_{14}=0$. Thus, $\rho_{ab}=$ for $a\neq b$. This
proves that ${\mathcal B}=0$.
 Therefore ${\mathcal R}=\displaystyle{\frac{s}{12}}Id$, i.e. $(M,g)$ is
of constant sectional curvature.

Conversely, if ${\mathcal R}=\displaystyle{\frac{s}{12}}Id$ a simple
straightforward computation shows that identity (\ref{Jklr}) is
satisfied. Hence ${\mathscr K}_1$ is integrable.

\end{proof}

\end{document}